\documentclass{article}
\usepackage{amsmath}
\usepackage{amsfonts,amssymb}
\input amssym.def

\numberwithin{equation}{section}

\newcommand{\End}{{\rm E n d }}
\newcommand{\Hom}{{\rm H o m }}
\newcommand{\Mat}{{\rm M a t }}

\newcommand{\ES}{{\rm E}{\cal S}}

\newcommand{\dia}{\diamond}

\newcommand{\tcM}{\widetilde{\cal M}}
\newcommand{\tcL}{\widetilde{\cal L}}
\newcommand{\tcK}{\widetilde{\cal K}}
\newcommand{\tcS}{\widetilde{\cal S}}

\newcommand{\tal}{\widetilde{\alpha}}
\newcommand{\tf}{\widetilde{f}}
\newcommand{\tg}{\widetilde{g}}
\newcommand{\tth}{\widetilde{h}}

\newcommand{\brarrow}{\succ\rightarrow}

\newcommand{\bbrarrow}{\succ\succ\rightarrow}
\newcommand{\bblarrow}{\leftarrow\prec\prec}

\newcommand{\SM}{{\cal S}M}

\newcommand{\Tp}{T_{poly}^{\bul}}

\newcommand{\cTp}{\cT^{\bul}_{poly}}

\newcommand{\cAb}{{\cal A}^{\bullet}}
\newcommand{\Omb}{\Om^{\bul}}

\newcommand{\OmT}{\Om^{\bul}(\cT^{\bullet}_{poly})}
\newcommand{\OmD}{\Om^{\bul}(\Cbu(\SM))}
\newcommand{\OmC}{\Om^{\bul}(\Cbd(\SM))}
\newcommand{\OmE}{\Om^{\bul}(\cE^{\bullet})}
\newcommand{\OM}{\cO_M}

\newcommand{\Cbu}{C^{\bullet}}
\newcommand{\Cbd}{C_{\bullet}}

\newcommand{\Linf}{L_{\infty}}

\newcommand{\tga}{\widetilde{\gamma}}

\newcommand{\al}{{\alpha}}
\newcommand{\la}{{\lambda}}
\newcommand{\h}{{\hbar}}
\newcommand{\bul}{{\bullet}}

\newcommand{\mF}{{\mathfrak{F}}}
\newcommand{\mG}{{\mathfrak{G}}}

\newcommand{\mb}{{\mathfrak{b}}}
\newcommand{\ms}{{\mathfrak{s}}}
\newcommand{\mg}{{\mathfrak{g}}}

\newcommand{\bE}{{\bar E}}

\newcommand{\Om}{{\Omega}}

\newcommand{\si}{{\sigma}}
\newcommand{\ga}{{\gamma}}

\newcommand{\G}{{\Gamma}}

\newcommand{\pa}{{\partial}}
\newcommand{\tr}{{\rm t r}}
\newcommand{\ch}{{\rm c h}}
\newcommand{\pr}{{\rm p r}}

\newcommand{\cK}{{\cal K}}
\newcommand{\cC}{{\cal C}}
\newcommand{\cM}{{\cal M}}

\newcommand{\cL}{{\cal L}}
\newcommand{\cD}{{\cal D}}
\newcommand{\cA}{{\cal A}}
\newcommand{\cG}{{\cal G}}

\newcommand{\cT}{{\cal T}}

\newcommand{\cS}{{\cal S}}
\newcommand{\cE}{{\cal E}}

\newcommand{\cO}{{\cal O}}

\newcommand{\bbC}{{\Bbb C}}
\newcommand{\bbR}{{\Bbb R}}
\newcommand{\bbZ}{{\Bbb Z}}

\newcommand{\La}{{\Lambda}}

\newcommand{\n}{{\nabla}}

\newcommand{\de}{{\delta}}

\newcommand{\tQ}{{\tilde{Q}}}

\newcommand{\trd}{{\rm t r d}}
\newcommand{\ttrd}{\widetilde{\rm t r d}}
\newcommand{\ind}{{\rm i n d}}
\newcommand{\cotr}{{\rm c o t r}}
\newcommand{\Lie}{{\rm L i e}}
\newcommand{\tind}{\widetilde{\rm i n d}}

\date{}
\newtheorem{defi}{Definition}

\newtheorem{lem}{Lemma}
\newtheorem{teo}{Theorem}

\newtheorem{pred}{Proposition}
\newtheorem{cond}{Condition}
\newtheorem{conj}{Conjecure}

\author{V.A. Dolgushev and V.N. Rubtsov}

\title{An algebraic index theorem
for Poisson manifolds}

\begin{document}

\large

\maketitle

\begin{abstract}
The formality theorem for Hochschild
chains of the algebra of functions on
a smooth manifold gives us a version of
the trace density map from the zeroth
Hochschild homology of a deformation quantization
algebra to the zeroth Poisson homology.
We propose a version of the algebraic index
theorem for a Poisson manifold which
is based on this trace density map.
\end{abstract}

~

\section{Introduction}
Various versions \cite{BNT}, \cite{BNT1}, \cite{Ch-I},
\cite{CFS}, \cite{Fedosov}, \cite{FFS}, \cite{NT} of
the algebraic index theorem
generalize the famous Atiyah-Singer index
theorem \cite{AS} from the case of a cotangent
bundle to an arbitrary symplectic manifold.
The first version of this theorem for Poisson
manifolds was proposed by D. Tamarkin and
B. Tsygan in \cite{TT}. Unfortunately, the proof of this
version is based on the formality conjecture for
cyclic chains \cite{Tsygan} which is not
yet established. In this paper we use the
formality theorem for Hochschild chains
\cite{FTHC}, \cite{Sh} to prove another
version of the algebraic index theorem
for an arbitrary Poisson manifold.
This version is based on the trace density
map from the zeroth Hochschild homology of
the deformation quantization algebra to
the zeroth Poisson homology of the manifold.

We denote by
$(M, \pi_1)$ a smooth real Poisson manifold and
by $\OM$ the algebra of smooth (real-valued)
functions on $M$. $TM$ (resp. $T^*M$) stands for
tangent (resp. cotangent) bundle of $M$\,.

Let $\OM^\h= (\OM[[\h]], *)$ be
a deformation quantization algebra of $(M,\pi_1)$
in the sense of \cite{Bayen} and \cite{Ber} and
let
\begin{equation}
\label{pi-h}
\pi = \h \pi_1 + \h^2 \pi_2 + \dots
\in \h\,\G(M, \wedge^2 TM)[[\h]]
\end{equation}
be a representative of Kontsevich's class
of $\OM^\h$\,.

One of the versions of the algebraic
index theorem for a symplectic manifold
\cite{Ch-I} describes a natural map
(see Eq. (31) and Theorem 4 in \cite{Ch-I})
$$
c l : K_0(\OM^\h) \to H^{top}_{DR}(M)((\h))
$$
from the $K$-theory of the deformation quantization
algebra $\OM^\h$ to the top degree De Rham
cohomology of $M$\,.
This map is obtained by composing
the trace density map \cite{FFS}
\begin{equation}
\label{trd}
\trd_{{\rm symp}} : \OM^\h/ [\OM^{\h}, \OM^{\h}]  \to
H^{top}_{DR}(M)((\h))
\end{equation}
from the zeroth Hochschild homology
$HH_0(\OM^\h) = \OM^\h/ [\OM^{\h}, \OM^{\h}]$
of $\OM^\h$ to
the top degree De Rham cohomology of $M$ with
the lowest component of the Chern character
(see example 8.3.6 in \cite{Loday})
\begin{equation}
\label{Ch}
\ch_{0,0} : K_0(\OM^\h) \to \OM^\h/ [\OM^{\h}, \OM^{\h}]
\end{equation}
from the $K$-theory of $\OM^\h$ to the
zeroth Hochschild homology of $\OM^\h$\,.

In the case of an arbitrary Poisson
manifold one cannot construct
the map (\ref{trd}). Instead, the formality
theorem for Hochschild chains \cite{FTHC}, \cite{Sh}
provides us with the map
\begin{equation}
\label{trd1}
\trd : \OM^\h/ [\OM^{\h}, \OM^{\h}]  \to
HP_0(M, \pi)\,,
\end{equation}
from the zeroth Hochschild homology of $\OM^\h$
to the zeroth Poisson homology
\cite{JLB}, \cite{Koszul} of $\pi$ (\ref{pi-h})\,.

According to J.-L. Brylinski \cite{JLB},
if $(M, \pi_1)$ is a
symplectic manifold
then we have
the following isomorphism:
$$
HP_{\bul}(M, \pi)[\h^{-1}] \cong
H^{dim\, M -\bul}_{DR}(M)((\h))\,.
$$
Thus, in the symplectic case
the map (\ref{trd}) can be obtained
from the map (\ref{trd1}).

For this reason we also refer to (\ref{trd1})
as {\it the trace density map}.

Composing (\ref{trd1}) with (\ref{Ch}) we
get the map
\begin{equation}
\label{ind}
\ind : K_0 (\OM^\h) \to HP_0(M, \pi) \,.
\end{equation}
Let us call this map {\it the quantum index
density}.

On the other hand setting $\h=0$ gives
us the obvious map
\begin{equation}
\label{prin}
\si :  K_0(\OM^{\h}) \to K_0(\OM)
\end{equation}
which we call {\it the principal symbol map}.

We recall that
\begin{pred}[J. Rosenberg, \cite{Ros}]
\label{ind-si}
The map $\ind$ (\ref{ind})
factors through the map $\si$ (\ref{prin}).
\end{pred}
The proof of this proposition is nice and transparent.
For this reason we decided to recall it here in the
introduction.

~\\
{\bf Proof.} First, recall that
for every associative algebra
$B$ a finitely generated projective module
can be represented by an idempotent in
the algebra $\Mat(B)$ of finite size matrices
with entries in $B$.

Next, let us show that the map
(\ref{prin}) is surjective. To do this,
it suffices to show that for every idempotent
$q$ in $\Mat(\OM)$ there exists an idempotent
$Q$ in $\Mat(\OM^\h)$ such that
$$
Q \Big|_{\h = 0} = q\,.
$$
The desired idempotent is produced by
the following equation
(see Eq. $(6.1.4)$ on page $185$ in
\cite{Fedosov}):
\begin{equation}
\label{to-samoe}
Q = \frac1{2} +
\Big( q  - \frac1{2}\Big) *
\Big(\, 1 + 4(q *q  - q)\, \Big)^{-1/2}\,.
\end{equation}
where the
last term in the right hand side
is understood as the expansion of the function
$y=x^{-1/2}$ in $4(q *q  - q)$ around
the point $(x=1, y=1)$\,. Since
$q$ is an idempotent in $\Mat(\OM)$
$$
q * q  - q  = 0 \quad mod \quad \h
$$
and therefore the expansion in (\ref{to-samoe})
makes sense.

A direct computation shows that the element
$Q$ defined by (\ref{to-samoe})
is indeed an idempotent in $\Mat(\OM^\h)$\,.

Thus, it suffices to show that if
two idempotents $P$ and $Q$ in $\Mat(\OM^{\h})$
have the same
principal part then $\ind([P])= \ind([Q])$\,.
Here $[P]$ (resp. $[Q]$) denotes the class
in $K_0(\OM^\h)$ represented by $P$ (resp. $Q$)

For this, we first show that if
\begin{equation}
\label{P-Q-h}
P\Big|_{\h=0} = Q \Big|_{\h=0}
\end{equation}
then $P$ can be connected to $Q$ by
a smooth path $P_t$ of idempotents in
$\Mat(\OM^\h)$\,.

Indeed if we define
the following smooth path
$$
P^0_t = ((1-t)P + t Q)
$$
in the algebra $\Mat(\OM^\h)$ and
plug $P^0_t$ into Equation
(\ref{to-samoe}) instead of $q$
we get the path of idempotents
in the algebra $\Mat(\OM^\h)$
\begin{equation}
\label{to-samoe1}
P_t = \frac1{2} +
\Big( P^0_t - \frac1{2}\Big) *
\Big(\, 1 + 4(P^0_t * P^0_t -
P^0_t )\, \Big)^{-1/2}\,,
\end{equation}
which connects $P$ with $Q$.
Due to Equation (\ref{P-Q-h}) the right hand
side of (\ref{to-samoe1}) is well defined
as a formal power series in $\h$\,.

Let us now show that $P$ and $Q$ represent
the same class in $K_0(\OM^\h)$. Since
$$
d_t P_t = (d_t P_t)* P_t
+ P_t *  ( d_t P_t)\,,
$$
$$
P_t * ( d_t P_t) * P_t = 0\,.
$$
Hence, for $P_t$ we have
$$
 d_t P(t) = [P(t), P(t) * (d_t P(t)) -
(d_t P(t)) * P(t)]\,.
$$
Therefore, since $P_t$ connects $P$ and $Q$\,,
$[P]= [Q]$ in $K_0(\OM^\h)$ and the
proposition follows. $\Box$

In this paper we propose a version of the
algebraic index theorem which describes how
the quantum index density (\ref{ind})
factors through the principal symbol map (\ref{prin}).

More precisely, using the generalization \cite{BGNT},
\cite{Endom}
of the formality theorem for Hochschild
chains \cite{FTHC}, \cite{Sh} to the algebra
of endomorphisms
of a vector bundle, we construct the map
\begin{equation}
\label{ind-c}
\ind_c : K_0(\OM) \to HP_0(M, \pi)\,,
\end{equation}
which makes the following diagram
\begin{equation}\label{ind-th}
\begin{array}{ccccc}
K_0 (\OM^\h) &~ &\stackrel{\ind}{\longrightarrow} & ~
& HP_0(M, \pi)\\
~ & \searrow^{\si} &  ~  &  ~^{\ind_c}\nearrow  & ~\\
 ~ & ~ &  K_0(\OM)  & ~ & ~
\end{array}
\end{equation}
commutative.

In this paper we refer
to the map $\ind_c$ as {\it the classical index density.}

The organization of the paper is as follows.
In the next section we fix notation and recall some results
we are going to use in this paper. In section 3 we prove some
useful facts about the twisting procedure of DGLAs and
DGLA modules by Maurer-Cartan elements. In section 4
we construct trace density map (\ref{trd1}), quantum
(\ref{ind}) and classical (\ref{ind-c}) index densities.
In section 5 we formulate and prove the main result of
this paper (see Theorem \ref{alg-ind}).
The concluding section consists of two parts.
In the first part
we describe the relation of our result to
the Tamarkin-Tsygan version \cite{TT} of
the algebraic index theorem. In the second
part we propose a
conjectural version of our index theorem
in the context of Rieffel's strict deformation
quantization \cite{Lands} of dual bundle of
a Lie algebroid.

~\\
{\bf Acknowledgment.}
We would like to thank P. Bressler, M. Kontsevich, B. Feigin,
G. Halbout, D. Tamarkin, and B. Tsygan for
useful discussions. We would like to thank the referee
for useful remarks and constructive
suggestions.
The results of this paper were
presented in the seminar on quantum groups and Poisson
geometry at Ecole Polytechnique.
We would like to thank
the participants of this seminar and especially
D. Sternheimer for questions and useful comments.
V.D. started this project when he was
a Liftoff Fellow of Clay Mathematics Institute
and he thanks this Institute for the
support. Bigger part of this paper was
written when V.D. was a Boas Assistant Professor
of Mathematics Department of Northwestern University.
V.D. thanks Northwestern University for perfect working
conditions and stimulating atmosphere.
During this project V.D. was a visitor of
FIB at ETH in Z\"urich and a visitor of the
University of Angers in France.
V.D. would like to thank
both institutes for hospitality and perfect
working conditions.
Both authors are partially supported by
the Grant for Support of Scientific
Schools NSh-8065.2006.2.
V.R. greatly acknowledge a partial
support of ANR-2005 (CNRS-IMU-RAS)
``GIMP'' and RFBR Grant N06-02-17382.
\section{Preliminaries}

In this section we fix
notation and recall some results we are going to
use in this paper.

For an associative algebra $B$
we denote by $\Mat_N(B)$
the algebra of $N\times N$ matrices over
$B$\,. The notation $\Cbd(B)$ is reserved for
the normalized Hochschild chain complex
of $B$ with coefficients in $B$
\begin{equation}
\label{chains}
\Cbd(B) = \Cbd(B,B)
\end{equation}
and the notation $\Cbu(B)$ is reserved for
the normalized Hochschild cochain complex of $B$
with coefficients in $B$ and with shifted
grading
\begin{equation}
\label{rule}
\Cbu(B) = C^{\bul + 1}(B,B)\,.
\end{equation}
The Hochschild coboundary operator
is denoted by $\pa$ and the
Hochschild boundary operator is
denoted by $\mb$\,.
We denote by $HH^{\bul}(B)$ the cohomology
of the complex $(\Cbu(B), \pa)$ and by
$HH_{\bul}(B)$ the homology of the complex
$(\Cbd(B), \mb)$\,.

It is well known that the
Hochschild cochain complex (\ref{rule})
carries the structure of a differential graded
Lie algebra. The corresponding
Lie bracket (see Eq. (3.2) on page 45 in \cite{thesis})
was originally introduced by M. Gerstenhaber
in \cite{Ger}. We will denote this bracket by
$[,]_G$\,.

The Hochschild chain complex (\ref{chains})
carries the structure of a differential graded
Lie algebra module over
the DGLA $\Cbu(B)$\,. We will denote the
action (see Eq. (3.5) on page 46 in \cite{thesis})
of cochains on chains by $R$\,.

The trace map $\tr$ \cite{Loday} is the map from the
Hochschild chain complex $\Cbd(\Mat_N(B))$ of
the algebra $\Mat_N(B)$ to the Hochschild chain
complex $\Cbd(B)$ of the algebra $B$\,. This map
is defined by the formula
\begin{equation}
\label{trace-map}
\tr (M_0 \otimes M_1 \otimes \dots \otimes M_k) =
\sum_{i_0, \dots, i_k}
(M_0)_{i_0 i_1} \otimes (M_1)_{i_1 i_2} \otimes
\dots \otimes (M_k)_{i_k i_0}\,,
\end{equation}
where $M_0, \dots, M_k$ are matrices in $\Mat_N(B)$
and $(M_a)_{ij}$ are the corresponding entries.

Dually, the cotrace map \cite{Loday}
$$
\cotr : \Cbu(B) \to \Cbu(\Mat_N(B))
$$
is defined by the formula
\begin{equation}
\label{cotrace-map}
(\cotr(P)(M_0, M_1, \dots , M_k))_{ij} =
\sum_{i_1, \dots, i_k}
P((M_0)_{i i_1}, (M_1)_{i_1 i_2},
\dots , (M_k)_{i_k  j})\,,
\end{equation}
where $P\in C^{k}(B)$ and
$M_0, \dots, M_k$ are, as above, matrices
in $\Mat_N(B)$\,.

``DGLA'' always means a differential graded Lie algebra.
The arrow $\brarrow$ denotes an $\Linf$-morphism
of DGLAs, the arrow $\bbrarrow$ denotes a
morphism of $\Linf$-modules,
and the notation
$$
\begin{array}{c}
\cL\\[0.3cm]
\downarrow_{\,mod}\\[0.3cm]
\cM
\end{array}
$$
means that $\cM$ is a DGLA module
over the DGLA $\cL$\,. The symbol
$\circ$ always stands for the composition of morphisms.
$\h$ denotes the formal deformation parameter.

Throughout this paper
$M$ is a smooth real ma\-ni\-fold.
For a smooth real vector bundle
$E$ over $M$, we denote by $\End(E)$
the algebra of endomorphisms of $E$\,.
For a sheaf $\cG$ of $\cO_M$-modules we denote by $\G(M, \cG)$
the vector space of global sections of $\cG$
and by $\Omb(\cG)$ the graded vector space
of exterior forms on $M$ with values in $\cG$\,.
In few cases, by abuse of notation, we
denote by $\Omb(\cG)$ the sheaf of exterior
forms with values in $\cG$\,. Similarly,
we sometimes refer to $\End(E)$ as the sheaf
of endomorphisms of a vector bundle $E$.
We specifically clarify the notation
when it is not clear from
the context.

$T^{\bul}_{poly}$ is the vector space of polyvector
fields with shifted grading
$$
T^{\bul}_{poly} = \G(M, \wedge^{\bul+1}_{\OM} TM)\,,
\qquad
T^{-1}_{poly} = \OM\,,
$$
and $\cAb$ is the graded vector
space of exterior forms:
$$
\cAb = \G(M,  \wedge^{\bul}_{\OM} T^*M)\,.
$$

$\Tp$ is the graded Lie algebra
with respect the so-called Schouten-Nijenhuis
bracket $[,]_{SN}$
(see Eq. (3.20) on page 50 in \cite{thesis}) and
$\cAb$ is the graded Lie algebra module
over $\Tp$ with respect to Lie derivative $L$
(see Eq. (3.21) on page 51
in \cite{thesis}). We will regard $\Tp$
(resp. $\cAb$) as the DGLA
(resp. the DGLA module)
with the zero differential.

Given a Poisson structure $\pi$
(\ref{pi-h}) one may introduce non-zero
differentials on the graded Lie algebra
$\Tp[[\h]]$ and on the graded Lie
algebra module $\cAb[[\h]]$\,. Namely,
$\Tp[[\h]]$ can be equipped with the Lichnerowicz
differential $[\pi, \, ]_{SN}$ \cite{Lich} and $\cAb[[\h]]$
can be equipped with the Koszul differential
$L_{\pi}$ \cite{Koszul}, where $L$ denotes the
Lie derivative. The cohomology of the
complex $(\Tp[[\h]], [\pi, \, ]_{SN} )$ is
called the Poisson cohomology of $\pi$.
For these cohomology groups
we reserve the notation $HP^{\bul}(M,\pi)$\,.
Similarly, the homology of the complex
$(\cAb[[\h]], L_{\pi})$ is called the Poisson
homology of $\pi$ and for the homology groups
of $(\cAb[[\h]], L_{\pi})$ we reserve notation
$HP_{\bul}(M, \pi)$\,.

We denote by $x^i$ local coordinates on $M$ and
by $y^i$ fiber coordinates in the tangent bundle
$TM$. Having these coordinates $y^i$ we
can introduce another local basis of exterior
forms $\{ dy^i \}$. We will use both bases
$ \{ dx^i \}$ and $\{ dy^i \}$. In particular,
the notation $\Omb(\cG)$
is reserved for the
$d y$-exterior forms with values in
the sheaf $\cG$ while $\cAb$ is reserved for
the $d x$-exterior forms.

We now briefly recall the
Fedosov resolutions (see Chapter 4 in \cite{thesis})
of polyvector fields, exterior
forms, and Hochschild (co)chains of $\cO_M$\,.
This construction has various
incarnations and it is referred to as the Gelfand-Fuchs
trick \cite{GF} or formal geometry \cite{G-Kazh}
in the sense of Gelfand and Kazhdan , or
mixed resolutions \cite{Ye} of Yekutieli.

We denote by
$\SM$ the formally completed symmetric algebra
of the cotangent bundle $T^*(M)$\,. Sections of the
sheaf $\SM$ can be viewed as formal power series in
tangent coordinates $y^i$\,. We regard $\SM$ as
the sheaf of algebras over $\OM$\,. In particular,
$\Cbu(\SM)$ is the sheaf of normalized
Hochschild cochains of $\SM$ over $\OM$. Namely,
the sections of $C^k(\SM)$ over an open subset
$U\subset M$ are $\OM$-linear polydifferential
operators with respect to the tangent
coordinates $y^i$
$$
P : \G(U,\SM)^{\otimes \, (k+1)} \to \G(U,\SM)
$$
satisfying the normalization condition
$$
P(\dots, f, \dots) = 0\,,
\qquad \forall ~ f \in \OM(U)\,.
$$
Similarly, $\Cbd(\SM)$ is the sheaf of normalized
Hochschild chains\footnote{In \cite{thesis} the sheaf
$\Cbu(\SM)$ is denoted by $\cD^{\bul}_{poly}$ and
the sheaf $\Cbd(\SM)$ is denoted by $\cC^{poly}_{\bul}$.}
of $\SM$ over $\OM$\,.
As in \cite{thesis} the tensor product in
$$
C_k(\SM) =
\underbrace{\SM \hat{\otimes}_{\OM} (\SM / \OM)
\hat{\otimes}_{\OM}
\dots
\hat{\otimes}_{\OM} (\SM / \OM)  }_{k+1}
$$
is completed in the adic topology in fiber
coordinates $y^i$ on the tangent bundle $TM$\,.

The cohomology of the complex of sheaves
$\Cbu(\SM)$ is the sheaf $\cTp$ of fiberwise
polyvector fields (see page 60 in \cite{thesis}).
The cohomology of the complex of sheaves
$\Cbd(\SM)$ is the sheaf $\cE^{\bul}$
of fiberwise differential
forms (see page 62 in \cite{thesis}). These are
$d x$-forms with values in $\SM$\,.

In \cite{thesis} (see Theorem 4 on page 68) it
is shown that the algebra $\Omb(\SM)$
can be equipped with a differential of the
following form
\begin{equation}
\label{DDD}
D = \n - \de + A\,,
\end{equation}
where
\begin{equation}
\label{nabla}
\n = dy^i \frac{\pa}{\pa x^i} -
dy^i \G^k_{ij}(x) y^j \frac{\pa}{\pa y^k}\,,
\end{equation}
is a torsion free connection with
Christoffel symbols $\G^k_{ij}(x)$,
\begin{equation}
\label{delta}
\de = dy^i \frac{\pa}{\pa y^i}\,,
\end{equation}
and
\begin{equation}
\label{A}
A=\sum_{p=2}^{\infty}dy^k A^j_{ki_1\dots i_p}(x)
y^{i_1} \dots
y^{i_p}\frac{\pa}{\pa y^j} \in \Om^1(\cT^0_{poly})\,.
\end{equation}
We refer to (\ref{DDD}) as the Fedosov differential.

Notice that $\de$ in (\ref{delta}) is also
a differential on
$\Omb(\SM)$ and (\ref{DDD}) can be viewed as a
deformation of $\de$ via the connection $\n$\,.

Let us recall from \cite{thesis}
the following operator
on\footnote{The arrow over $\pa$ in (\ref{del-1})
means that we use the left derivative with respect to
the ``anti-com\-muting'' variable $d y^k$.} $\Omb(\SM)$
\begin{equation}
\de^{-1}(a) =
\begin{cases}
\begin{array}{cc}
\displaystyle
y^k \frac {\vec{\partial}} {\partial (d y^k)}
\int\limits_0^1 a(x,t y,t d y)\frac{d t} t, & {\rm if}~
a \in \Om^{>0}(\SM)\,,\\
0, & {\rm otherwise}\,.
\end{array}
\end{cases}
\label{del-1}
\end{equation}
This operator satisfies the following
properties:
\begin{equation}
\de^{-1} \circ  \de^{-1} =0
\label{nilp}
\end{equation}
\begin{equation}
a=\chi(a) +\delta \delta ^{-1}a + \delta ^{-1}\delta a\,,
\qquad
\forall ~ a\in \Omb(\SM)
\label{Hodge}
\end{equation}
where
\begin{equation}
\label{chi}
\chi (a)= a\Big|_{y^i=dy^i=0}\,.
\end{equation}
These properties are used in
the proof of the
acyclicity of $\de$ and $D$ in positive
dimension.

According to Proposition $10$ on
page $64$ in \cite{thesis} the sheaves $\cTp$, $\Cbu(\SM)$,
$\cE^{\bul}$, and $\Cbd(\SM)$ are equipped with the
canonical action of the sheaf of Lie algebras
$\cT_{poly}^0$ and this action is compatible
with the corresponding (DG) algebraic structures.
Using this action in
chapter $4$ of \cite{thesis} we extend
the Fedosov differential (\ref{DDD}) to
a differential on the DGLAs (resp. DGLA modules)
$\Omb(\cTp)$, $\Omb(\cE^{\bul})$,
$\Omb(\Cbd(\SM))$, and $\Omb(\Cbu(\SM))$\,.

Using acyclicity of the Fedosov differential (\ref{DDD})
in positive dimension
one constructs in \cite{thesis} embeddings
of DGLAs and DGLA
modules\footnote{See Eq. (5.1) on page 81 in \cite{thesis}.}
\begin{equation}
\begin{array}{ccc}
\Tp &\stackrel{\la_T}{\,\longrightarrow\,} &
(\OmT, D, [,]_{SN})\\[0.3cm]
\downarrow^{L}_{\,mod}  & ~  &
\downarrow^{L}_{\,mod}\\[0.3cm]
\cAb  &\stackrel{\la_{\cA}}{\,\longrightarrow\,} &
(\OmE, D),
\end{array}
\label{diag-T}
\end{equation}

\begin{equation}
\begin{array}{ccc}
(\OmD, D+\pa, [,]_{G}) &\stackrel{\,\la_D}{\,\longleftarrow\,}
& \Cbu(\OM)\\[0.3cm]
\downarrow^{R}_{\,mod}  & ~  &     \downarrow^{R}_{\,mod} \\[0.3cm]
(\OmC, D+\mb) &\stackrel{\,\la_C}{\,\longleftarrow\,} &  \Cbd(\OM),
\end{array}
\label{diag-D}
\end{equation}
and shows that these are quasi-isomorphisms of the
corresponding complexes.

Furthermore, using Kontsevich's and Shoikhet's
formality theorems for $\bbR^d$ \cite{K}, \cite{Sh}
in \cite{thesis} one constructs the following
diagram
\begin{equation}
\label{diag-K-Sh}
\begin{array}{ccc}
(\OmT, D, [,]_{SN}) & \stackrel{\cK}{\brarrow} &
(\OmD, D+\pa, [,]_{G}) \\[0.3cm]
\downarrow^{L}_{\,mod}  & ~  &
\downarrow^{R}_{\,mod} \\[0.3cm]
(\OmE, D) & \stackrel{\cS}{\bblarrow} &
(\OmC, D+\mb)
\end{array}
\end{equation}
where $\cK$ is an $\Linf$ quasi-isomorphism
of DGLAs, $\cS$ is a quasi-isomorphism
of $\Linf$-modules  over
the DGLA $(\OmT, D, [,]_{SN})$\,,
and the $\Linf$-module structure on $\OmC$ is obtained
by composing the $\Linf$ quasi-isomorphism
$\cK$ with the DGLA modules structure $R$
(see Eq. (3.5) on p. 46 in \cite{thesis} for the
definition of $R$)\,.

Diagrams (\ref{diag-T}), (\ref{diag-D}) and (\ref{diag-K-Sh})
show that the DGLA module $\Cbd(\OM)$ of Hochschild
chains of $\OM$ is quasi-isomorphic to
the graded Lie algebra module $\cAb$ of its
cohomology.

~\\
{\bf Remark 1.} As in \cite{thesis} we use adapted versions
of Hochschild (co)chains for the algebras $\OM$ and $\End(E)$
of functions and of endomorphisms of a vector bundle $E$,
respectively. Thus, $\Cbu(\OM)$ is the complex of
 polydifferential
operators (see page 48 in \cite{thesis}) satisfying
the corresponding normalization condition.
$\Cbu(\End(E))$ is the complex of (normalized)
polydifferential operators acting on $\End(E)$
with coefficients in $\End(E)$\,.
Furthermore,  $\Cbd(\OM)$ is the complex of
(normalized) polyjets
$$
C_k(\OM) = \Hom_{\OM} (C^{k-1}(\OM), \OM)\,,
$$
and
$$
C_k(\End(E)) = \Hom_{\End(E)} (C^{k-1}(\End(E)), \End(E))\,.
$$
We have to warn the reader that the
complex $\Cbu(\OM)$
(resp. $\Cbu(\End(E))$) does not coincide with
the complex of Hochschild cochains of the algebra
$\OM$ of functions (resp. the algebra
$\End(E)$ of endomorphisms of $E$).
Similar expectation is wrong for Hochschild
chains. Instead we have the following inclusions:
$$
\Cbu(\OM) \subset C^{\bul}_{{\rm genuine}}(\OM)\,,
\qquad
\Cbu(\End(E)) \subset C^{\bul}_{{\rm genuine}}(\End(E))\,,
$$
$$
C_{\bul}^{{\rm genuine}}(\OM) \subset
\Cbd(\OM)\,,
\qquad
C_{\bul}^{{\rm genuine}}(\End(E)) \subset
\Cbd(\End(E)) \,,
$$
where $C_{\bul}^{{\rm genuine}}$ and
$C^{\bul}_{{\rm genuine}}$
refer to the original definitions (\ref{chains}),
(\ref{rule})
of Hochschild
chains and cochains, respectively.

~\\
{\bf Remark 2.} Unlike in \cite{thesis} we use
only normalized Hochschild (co)chains. It is not hard
to check that the results we need from \cite{thesis},
\cite{K}, and \cite{Sh}
also hold when this normalization condition is imposed.

Let us now recall the construction of
paper \cite{Endom}, in which the formality
of the DGLA module $(\,\Cbu(\End(E)), \Cbd(\End(E))\,)$
is proved.

The construction of \cite{Endom} is based
on the use of the following auxiliary sheaf
of algebras:
\begin{equation}
\label{ES}
\ES = \End(E) \otimes_{\OM} \SM
\end{equation}
considered as a sheaf of algebras over $\OM$.

It is shown in \cite{Endom} that the Fedosov
differential (\ref{DDD}) can be extended to
the following differential on $\Omb(\ES)$
\begin{equation}
\label{DDD-E1}
D^E = D + [\ga^E,\, ]\,, \qquad
\ga^E= \G^E + \tga^E\,,
\end{equation}
where $\G^E$ is a connection form of $E$ and
$\tga^E $ is an element in $\Om^1(\ES)$
defined by iterating the equation:
\begin{equation}
\label{iter-ga-E}
\ga^E = \G^E + \de^{-1} (\n \ga^E
+ A (\ga^E) + \frac{1}{2} [\ga^E, \ga^E])
\end{equation}
in degrees in fiber coordinates $y^i$ of
the tangent bundle $TM$\,.

The differential (\ref{DDD-E1}) naturally
extends to the DGLA $\Omb(\Cbu(\ES))$ and
to the DGLA module $\Omb(\Cbd(\ES))$\,. Namely,
on $\Omb(\Cbu(\ES))$ the differential $D^E$
is defined by the formula
\begin{equation}
\label{DDD-E-coch}
D^E = D + [\pa \ga^E,\, ]_G\,,
\end{equation}
and on $\Omb(\Cbd(\ES))$ it is defined by
the equation
\begin{equation}
\label{DDD-E-ch}
D^E = D + R_{\pa \ga^E }\,.
\end{equation}
Here, $\pa$ is the Hochschild coboundary
operator and $R$ denotes the actions of
cochains on chains.

Then,  generalizing the construction
of the maps $\la_D$ and $\la_C$ in
(\ref{diag-D}) one gets the
following embeddings
of DGLAs and their modules
\begin{equation}
\begin{array}{ccc}
(\Omb(\Cbu(\ES)), D^E + \pa, [,]_{G}) &
\stackrel{\,\la^E_D}{\,\longleftarrow\,} & \Cbu(\End(E)) \\[0.3cm]
\downarrow^{R}_{\,mod}  & ~  &     \downarrow^{R}_{\,mod} \\[0.3cm]
(\Omb(\Cbd(\ES)), D^E + \mb) &\stackrel{\,\la^E_C}{\,\longleftarrow\,}
&   \Cbd(\End(E))\,.
\end{array}
\label{diag-C-E}
\end{equation}
Similarly to Propositions 7, 13, 15 in \cite{thesis}
one can easily show that $\la^E_D$ and $\la^E_C$ are
quasi-isomorphisms of the corresponding complexes.

Finally, the DGLA modules
$(\,\Omb(\Cbu(\ES)), \Omb(\Cbd(\ES))\,)$ and
$(\,\Omb(\Cbu(\SM)), \Omb(\Cbd(\SM))\,)$ are
connected in \cite{Endom} by the following
commutative diagram
of quasi-isomorphisms of DGLAs and their modules
\begin{equation}
\begin{array}{ccc}
(\Omb (\Cbu(\SM)), D + \pa, [,]_{G})
&
\stackrel{\cotr^{tw}}{\,\rightarrow\,}
&
(\Omb (\Cbu(\ES)), D^E + \pa, [,]_{G}) \\[0.3cm]
 \downarrow_{\,mod}
&
~
&
 \downarrow_{\,mod} \\[0.3cm]
(\Omb (\Cbd(\SM)), D+\mb) &
\stackrel{\tr^{tw}}{\,\leftarrow\,}
&
(\Omb (\Cbd(\ES)), D^E + \mb)\,,
\end{array}
\label{tr-cotr-tw}
\end{equation}
where
\begin{equation}
\label{oni1}
\cotr^{tw} = \exp (-[\ga^E\,, \, ]_G) \circ \cotr\,,
\qquad
\tr^{tw} = \tr \circ \exp (R_{\ga^E})\,,
\end{equation}
and $\tr$, $\cotr$ are the maps
\begin{equation}
\label{tr-cotr}
\tr : \Cbd(\ES) \to \Cbd(\SM)\,,
\qquad
\cotr : \Cbu(\SM) \to \Cbu(\ES)\,
\end{equation}
defined as in (\ref{trace-map}) and
(\ref{cotrace-map}).

It should be remarked that the element
$\ga^E$ in (\ref{oni1}) comprises the connection
form of $E$ (\ref{DDD-E1}) and hence
can be viewed as a section of the
sheaf $\Om^1(C^{-1}(\ES))$ only locally.
The compositions $\tr^{tw}$ and
$\cotr^{tw}$ are still well defined due to the fact
that we consider normalized Hochschild
(co)chains.

Diagrams (\ref{diag-T}), (\ref{diag-K-Sh}),
(\ref{diag-C-E}), and (\ref{tr-cotr-tw})
give us the desired chain
of formality quasi-isomorphisms for the
DGLA module  $(\,\Cbu(\End(E)), \Cbd(\End(E))\,)$\,.

\section{The twisting procedure revisited}
In this section we will prove some
general facts about the twisting
by a Maurer-Cartan element.
See Section $2.4$ in \cite{thesis} in
which this procedure is discussed in more
details.

Let $(\cL, d_{\cL} , [\,,\,]_{\cL})$
$(\tcL, d_{\tcL} , [\,,\,]_{\tcL})$
and be two DGLAs over $\bbR$.
Since we deal with
deformation theory questions it is convenient
for our purposes to extend the field $\bbR$
to the ring $\bbR[[\h]]$ from the very beginning
and consider $\bbR[[\h]]$-modules
$\cL[[\h]]$ and $\tcL[[\h]]$
with DGLA structures extended
from $\cL$ and $\tcL$ in the obvious way.

By definition $\al$ is a Maurer-Cartan element
of the DGLA $\cL[[\h]]$ if $\al\in \h\,\cL[[\h]]$\,,
it has degree $1$ and satisfies
the equation:
\begin{equation}
\label{MC-def}
d_{\cL} \al + \frac1{2} [\al, \al]_{\cL} = 0\,.
\end{equation}
The first two conditions can be written concisely
as $\al\in \h\, \cL^1[[\h]]$\,.

Notice that, $\mg(\cL) = \h\, \cL^0[[\h]]$ forms an
ordinary (not graded) Lie algebra over $\bbR[[\h]]$\,.
Furthermore, $\mg(\cL)$ is obviously pronilpotent and
hence can be exponentiated to the group
\begin{equation}
\label{mG}
\mG(\cL) = \exp (\, \h\, \cL^0[[\h]] \, )\,.
\end{equation}
The natural action of this group on $\cL[[\h]]$
can be introduced by exponentiating the adjoint
action of $\mg(\cL)$\,.

The action of $\mG(\cL)$ on the  Maurer-Cartan elements
of the DGLA $\cL[[\h]]$ is given by the formula
\begin{equation}
\label{action}
\exp(\xi)[\al] = \al + f([\,\,,\xi]_{\cL}) (d_{\cL}\xi +
[\al, \xi]_{\cL})\,,
\end{equation}
where $\al$ is a Maurer-Cartan element of $\cL[[\h]]$\,,
$\xi \in \h \cL^0[[\h]]$\,, $f(x)$ is the function
$$
f(x) = \frac{e^x-1}{x}
$$
and the expression $f([\,\,,\xi]_{\cL})$ is defined
via the Taylor expansion of $f(x)$ around the point
$x=0$\,.

We call Maurer-Cartan elements equivalent
if they lie on
the same orbit of the action (\ref{action}).
The set of these orbits
is called the moduli space of the
Maurer-Cartan elements.

We have the following proposition:
\begin{pred}[W. Goldman and J. Millson, \cite{Gold}]
\label{Goldman}
If $f$ is a quasi-isomorphism from the DGLA
$\cL$ to the DGLA $\tcL$ then the induced map
between the moduli spaces of Maurer-Cartan elements of
$\cL[[\h]]$ and $\tcL[[\h]]$ is an isomorphism
of sets.
\end{pred}

Every Maurer-Cartan element $\al$ of $\cL[[\h]]$ can be
used to modify the DGLA structure
on $\cL[[\h]]$\,. This modified structure is called \cite{Q}
the DGLA structure twisted by the Maurer-Cartan
$\al$\,. The Lie bracket of the twisted DGLA structure
is the same and the differential is given by
the formula:
\begin{equation}
\label{d-al}
d^{\al}_{\cL} = d_{\cL} + [\al, \, ]_{\cL}\,.
\end{equation}
We will denote the DGLA $\cL[[\h]]$ with the bracket
$[\,,\,]_{\cL}$ and the differential $d_{\cL}^{\al}$ by
$\cL^{\al}$\,.

Two DGLAs are called quasi-isomorphic if there
is a chain of quasi-isomorphisms
$f$, $f_1$, $f_2$, $\dots$, $f_n$ connecting
$\cL$ with $\tcL$:
\begin{equation}
\label{cL-tcL}
\cL \,\stackrel{f}{\rightarrow}\, \cL_1
\,\stackrel{f_1}{\leftarrow}\,
\cL_2   \,\stackrel{f_2}{\rightarrow}\,
\dots
 \,\stackrel{f_{n-1}}{\rightarrow}\,
\cL_n \,\stackrel{f_n}{\leftarrow}\,  \tcL
\end{equation}

This chain naturally extends to the chain of
quasi-isomorphisms of DGLAs over $\bbR[[\h]]$
\begin{equation}
\label{cL-tcL-h}
\cL[[\h]] \,\stackrel{f}{\rightarrow}\, \cL_1[[\h]]
\,\stackrel{f_1}{\leftarrow}\,
\cL_2[[\h]]   \,\stackrel{f_2}{\rightarrow}\,
\dots
 \,\stackrel{f_{n-1}}{\rightarrow}\,
\cL_n[[\h]] \,\stackrel{f_n}{\leftarrow}\,  \tcL[[\h]]
\end{equation}

We would like to prove that
\begin{pred}
\label{chain-twist}
For every Maurer-Cartan element $\al$ of
$\cL[[\h]]$ the chain of quasi-isomorphisms (\ref{cL-tcL-h})
can be upgraded to the chain
\begin{equation}
\label{cL-tcL-tw}
\cL^{\al} \,\stackrel{f}{\rightarrow}\, \cL^{\al_1}_1
\,\stackrel{\tf_1}{\leftarrow}\,
\cL^{\al_2}_2   \,\stackrel{\tf_2}{\rightarrow}\,
\dots
 \,\stackrel{\tf_{n-1}}{\rightarrow}\,
\cL^{\al_n}_n \,\stackrel{\tf_n}{\leftarrow}\,  \tcL^{\tal}\,,
\end{equation}
where $\al_i$ (resp. $\tal$) are Maurer-Cartan
elements of $\cL_i[[\h]]$ (resp. $\tcL[[\h]]$)
and the quasi-isomorphisms $\tf_i$ are obtained
from $f_i$ by composing with the action of
an element in the group $\mG(\cL_i)$\,.
\end{pred}
{\bf Proof} runs by induction on $n$\,.
We, first, prove the base of the
induction $(n=1)$ and then the step
follows easily from the statement of the
proposition for $n=1$\,.

We set $\al_1$ to be
$$
\al_1 = f(\al)\,.
$$
Since $f_1$ is a quasi-isomorphism from $\tcL$
to $\cL_1$, by Proposition \ref{Goldman}, there exists
a Maurer-Cartan element $\tal \in \tcL$ such that
$f_1(\tal)$ is equivalent to $\al_1$\,.

Let $T_1$ be an element of the group $\mG(\cL_1)$
which transforms $f_1(\tal)$ to $\al_1$\,.
Thus by setting
\begin{equation}
\label{tf-1}
\tf_1 =  T_1 \circ f_1
\end{equation}
we get the desired chain for $n=1$
\begin{equation}
\label{cL-tcL-tw1}
\cL^{\al} \,\stackrel{f}{\rightarrow}\, \cL^{\al_1}_1
\,\stackrel{\tf_1}{\leftarrow}\, \tcL^{\tal}
\end{equation}
and the proposition follows. $\Box$

Chain (\ref{cL-tcL-tw}) gives us an isomorphism
\begin{equation}
\label{I-tal}
I_{\tal} : H^{\bul}(\tcL^{\tal})
\stackrel{\sim}{\rightarrow}  H^{\bul}(\cL^{\al})
\end{equation}
from the cohomology of the DGLA $\tcL^{\tal}$
to the cohomology of the DGLA $\cL^{\al}$\,.
This isomorphism depends on choices of
Maurer-Cartan elements in the intermediate terms $\cL_i[[\h]]$
and the choices of elements from the groups
$\mG(\cL_i)$\,.

We claim that
\begin{pred}
\label{abs-nonsense}
If $\cL$ and $\tcL$ are DGLAs connected by
the chain of quasi-isomorphisms (\ref{cL-tcL}),
$\al$ is a Maurer-Cartan element in
$\h \cL[[\h]]$ and
\begin{equation}
\label{cL-tcL-tw11}
\cL^{\al} \,\stackrel{f}{\rightarrow}\, \cL^{\al'_1}_1
\,\stackrel{\tf'_1}{\leftarrow}\,
\cL^{\al'_2}_2   \,\stackrel{\tf'_2}{\rightarrow}\,
\dots
 \,\stackrel{\tf'_{n-1}}{\rightarrow}\,
\cL^{\al'_n}_n \,\stackrel{\tf'_n}{\leftarrow}\,  \tcL^{\tal'}\,,
\end{equation}
is another chain of quasi-isomorphism obtained
according to Proposition \ref{chain-twist} then
there exists and element $T_{\tcL}$ of the group
$\mG(\tcL)$ such that
\begin{equation}
\label{T-tcL}
T_{\tcL}(\tal) = \tal' \,,
\end{equation}
and
\begin{equation}
\label{I-T-tcL}
I_{\tal} = I_{\tal'} \circ T_{\tcL} \,.
\end{equation}
\end{pred}
Before proving the proposition let us consider the
case $n=0$ with a slight modification.
More precisely, we consider a quasi-isomorphism
$g= f_0$ from the DGLA $\tcL$ to the DGLA $\cL$
\begin{equation}
\label{g}
\cL \,\stackrel{g}{\leftarrow}\, \tcL
\end{equation}
and suppose that $\al$ and $\al'$ are two equivalent
Maurer-Cartan elements of $\cL[[\h]]$\,.

Due to Proposition \ref{Goldman} there exist
Maurer-Cartan elements $\tal$ and $\tal'$ in
$\tcL[[\h]]$ and elements $T$ and $T'$ of the group
$\mG(\cL)$ such that
$$
T(g(\tal)) = \al\,, \qquad
T'(g(\tal')) = \al'\,.
$$
Hence, by setting
\begin{equation}
\label{tg-tg'}
\tg = T \circ g\,, \qquad
\tg' = T' \circ g
\end{equation}
we get the following pair of quasi-isomorphisms
of twisted DGLAs
\begin{equation}
\label{tg}
\cL^{\al} \, \stackrel{\tg}{\leftarrow} \, \tcL^{\tal}\,,
\end{equation}
\begin{equation}
\label{tg'}
\cL^{\al'} \, \stackrel{\tg'}{\leftarrow} \, \tcL^{\tal'}\,.
\end{equation}
Let us prove the following auxiliary statement:

\begin{lem}
\label{n0}
If $\al$ and $\al'$ are Maurer-Cartan elements
of $\cL[[\h]]$ connected by the action of
an element $T_{\cL}\in \mG(\cL)$
$$
T_{\cL} (\al) = \al'\,,
$$
then there exists
an element $T_{\tcL}$ of the group $\mG(\tcL)$
such that
$$
T_{\tcL} (\tal) =  \tal'
$$
and the diagram of DGLAs
\begin{equation}
\label{diag-T-tcL}
\begin{array}{ccc}
\cL^{\al} & \,\stackrel{\tg}{\leftarrow}\, & \tcL^{\tal} \\[0.3cm]
 \downarrow^{T_{\cL}} &  ~  &  \downarrow^{T_{\tcL}}     \\[0.3cm]
\cL^{\al'} & \,\stackrel{\tg'}{\leftarrow}\, & \tcL^{\tal'}
\end{array}
\end{equation}
commutes up to homotopy.
\end{lem}
{\bf Proof.} Since the Maurer-Cartan element $g(\tal)$ and
$g(\tal')$ are equivalent, then so are the
Maurer-Cartan elements $\tal$ and $\tal'$\,.
Hence, there exists an element $T_0\in \mG(\tcL)$ such that
$$
T_0(\tal) = \tal'\,.
$$

The following diagram shows the relations between
various Maurer-Cartan elements in question:
\begin{equation}
\label{T-alpha}
\begin{array}{ccccc}
\al & \stackrel{T}{\leftarrow} & g(\tal) &
\stackrel{g}{\leftarrow} & \tal \\[0.3cm]
\downarrow^{T_{\cL}} & ~ & \downarrow^{g(T_0)}
& ~ &\downarrow^{T_0} \\[0.3cm]
\al' & \stackrel{T'}{\leftarrow} & g(\tal') &
\stackrel{g}{\leftarrow} & \tal'\,.
\end{array}
\end{equation}
From this diagram we see that
$$
T_{\cL} T [g(\tal)] = T' g(T_0) [g(\tal)]\,,
$$
or equivalently
\begin{equation}
\label{TTTT}
g(T_0)^{-1} (T')^{-1} T_{\cL} T [g(\tal)] = [g(\tal)]\,.
\end{equation}
Therefore the element $g(T_0)^{-1} (T')^{-1} T_{\cL} T $
belongs to the subgroup $\mG(\cL,  g(\tal))\subset \mG(\cL)$
of elements preserving $g(\tal)$\,.

It is not hard to show that
the subgroup $\mG(\cL,  g(\tal))$ is
\begin{equation}
\label{MC-fix1}
\mG(\cL,  g(\tal)) = \exp (\h \cL^0 [[\h]] \cap
\ker (d_{\cL} + [g(\tal), \, ]_{\cL})\,.
\end{equation}
In other words, there exists a
$d_{\cL} + [g(\tal), \, ]_{\cL}$-closed
element $\xi \in \h \cL^0[[\h]]$ such that
$$
g(T_0)^{-1} (T')^{-1} T_{\cL} T  = \exp (\xi)\,.
$$

Since the map $g$ from $\tcL^{\tal}$ to
$\cL^{g(\tal)}$ is
a quasi-isomorphism, there exists a element
$\psi\in \h \tcL^0[[\h]]$ such that
$$
d_{\tcL}\psi + [\tal, \psi]_{\tcL} = 0\,,
$$
and the difference $\xi - g(\psi)$ is
$d_{\cL} + [g(\tal),\, ]_{\cL}$ - exact.
Therefore the elements
$$
T_{\cL} T
$$
and
$$
T' g(T_0\, \exp(\psi))
$$
induce homotopic maps from the DGLA
$\cL^{g(\tal)}$ to the DGLA $\cL^{\al'}$\,.

Thus, $T_{\tcL} = T_0 \exp(\psi)$ is the desired
element of the group $\mG(\tcL)$ which makes
Diagram (\ref{diag-T-tcL}) commutative up to homotopy.

The lemma is proved. $\Box$

~\\
{\bf Proof of Proposition \ref{abs-nonsense}.}
The proof runs by induction on $n$ and the
base of the induction $n=0$ follows immediately
from Lemma \ref{n0}.

If the statement is proved for $n=2k$ then the
statement for $n=2k+1$ also follows from Lemma
\ref{n0}.

Since the source of the map $f_{2k}$ is $\cL_{2k}$
the case $n=2k$ follows immediately from the
case $n=2k-1$ and this concludes the
proof of the proposition. $\Box$

For a DGLA module $\cM$ over a DGLA $\cL$
the direct sum $\cL \oplus \cM$ carries the natural
structure of a DGLA. Namely, this DGLA is the
semi-direct product of $\cL$ and $\cM$ where
$\cM$ is viewed as a DGLA with the zero
bracket. Morphisms between two such semi-direct products
$\cL \oplus \cM $ and $\tcL \oplus \tcM$ correspond to
morphisms between the DGLA modules $(\cL,\cM)$ and
$(\tcL, \tcM)$\,.  Furthermore,
 twisting the DGLA structure on
$\cL[[\h]] \oplus \cM[[\h]]$ by a Maurer-Cartan
element $\al\in \h\, \cL[[\h]]$ gives us the semi-direct
product $\cL^{\al} \oplus \cM^{\al}$ corresponding to the
DGLA module $\cM^{\al}$ with the differential twisted
by the action of the Maurer-Cartan element $\al$\,.

This observation allows us to generalize
Propositions \ref{chain-twist} and \ref{abs-nonsense}
to a chain of quasi-iso\-mor\-phisms of DGLA modules:
\begin{equation}
\label{chain-mod}
(\cL, \cM) \,\stackrel{h}{\rightarrow}\,
(\cL_1, \cM_1) \,\stackrel{h_1}{\leftarrow}\,
(\cL_2, \cM_2) \,\stackrel{h_2}{\rightarrow}\,
\dots
 \,\stackrel{h_{n-1}}{\rightarrow}\,
(\cL_n, \cM_n) \,\stackrel{h_n}{\leftarrow}\,
(\tcL, \tcM)\,.
\end{equation}
Namely,
\begin{pred}
\label{chain-mod-twist}
For every Maurer-Cartan element $\al$ of
$\cL[[\h]]$ the chain of quasi-isomorphisms (\ref{chain-mod})
can be upgraded to the chain
\begin{equation}
\label{cL-cM-tw}
(\cL^{\al}, \cM^{\al} ) \,\stackrel{h}{\rightarrow}\,
(\cL^{\al_1}_1, \cM^{\al_1}_1)
\,\stackrel{\tth_1}{\leftarrow}\,
(\cL^{\al_2}_2, \cM^{\al_2}_2)   \,\stackrel{\tth_2}{\rightarrow}\,
\dots
 \,\stackrel{\tth_{n-1}}{\rightarrow}\,
(\cL^{\al_n}_n, \cM^{\al_n}_n)
\,\stackrel{\tth_n}{\leftarrow}\,
(\tcL^{\tal}, \tcM^{\tal})\,,
\end{equation}
where $\al_i$ (resp. $\tal$) are Maurer-Cartan
elements of $\cL_i[[\h]]$ (resp. $\tcL[[\h]]$)
and the quasi-isomorphisms $\tth_i$ are obtained
from $h_i$ by composing with the action of
an element in the group $\mG(\cL_i)$\,. $\Box$
\end{pred}
Chain (\ref{cL-cM-tw}) gives us an isomorphism
\begin{equation}
\label{J-tal}
J_{\tal} : H^{\bul}(\tcM^{\tal})
\stackrel{\sim}{\rightarrow}  H^{\bul}(\cM^{\al})
\end{equation}
from the cohomology of the DGLA module $\tcM^{\tal}$
to the cohomology of the DGLA module $\cM^{\al}$\,.
This isomorphism depends on choices of
Maurer-Cartan elements in the DGLAs $\cL_i[[\h]]$
and the choices of elements from the groups
$\mG(\cL_i)$\,.

We claim that
\begin{pred}
\label{abs-nonsense-mod}
If $(\cL, \cM)$ and $(\tcL, \tcM)$ are DGLA
modules connected by
the chain of quasi-isomorphisms (\ref{chain-mod}),
$\al$ is a Maurer-Cartan element in
$\h \cL[[\h]]$ and
\begin{equation}
\label{cL-cM-tw11}
(\cL^{\al}, \cM^{\al}) \,\stackrel{h}{\rightarrow}\,
(\cL^{\al'_1}_1, \cM^{\al'_1}_1)
\,\stackrel{\tth'_1}{\leftarrow}\,
(\cL^{\al'_2}_2, \cM^{\al'_2}_2)
\,\stackrel{\tth'_2}{\rightarrow}\,
\dots
 \,\stackrel{\tth'_{n-1}}{\rightarrow}\,
(\cL^{\al'_n}_n, \cM^{\al'_n}_n)
\,\stackrel{\tth'_n}{\leftarrow}\,
(\tcL^{\tal'}, \tcM^{\tal'}) \,,
\end{equation}
is another chain of quasi-isomorphism obtained
according to Proposition \ref{chain-mod-twist} then
there exists an element $T_{\tcL}$ of the group
$\mG(\tcL)$ such that
\begin{equation}
\label{T-tcL-mod}
T_{\tcL}(\tal) = \tal' \,,
\end{equation}
and
\begin{equation}
\label{J-T-tcL}
J_{\tal} = J_{\tal'} \circ T_{\tcL} \,.  \qquad \Box
\end{equation}
\end{pred}

\section{Trace density, quantum and classical
index densities}

In this section we recall the construction
of the trace density map (\ref{trd1}) which gives
the quantum index density (\ref{ind}).
We also construct the classical index
density (\ref{ind-c}) using the chain
(\ref{diag-T}), (\ref{diag-K-Sh}),
(\ref{diag-C-E}), (\ref{tr-cotr-tw})
of formality quasi-isomorphisms for $\Cbd(\End(E))$\,.

Let, as above, $\OM^{\h}= (\OM[[\h]], *)$ be
a deformation quantization algebra of
the Poisson manifold $(M, \pi_1)$ and
$\pi$ (\ref{pi-h}) be a representative of
Kontsevich's class of $\OM^\h$\,.

It can be shown that every star-product
$*$ is equivalent to the so-called
natural star-product \cite{GR}. These
are the star-products
$$
a * b = a\, b +
\sum_{k=1}^{\infty} \h^k B_k(a,b)
$$
for which the bidifferential operators $B_k$
satisfy the following condition:
\begin{cond}
\label{uslovie}
For all $k\ge 1$ the bidifferential
operator $B_k$ has the
order at most $k$ in each argument.
\end{cond}
In this paper we tacitly assume that
the above condition holds for the
star-product $*$\,.

Using the map $\la_T$ from (\ref{diag-T})
we lift $\pi$ to the Maurer-Cartan element
$\la_T(\pi)$ in the DGLA $\Omb(\cTp)[[\h]]$\,.
This element allows us to extend the differential
$D$ (\ref{DDD}) on $\OmT[[\h]]$ to
\begin{equation}
\label{la-pi}
D + [\la_T(\pi),\,]_{SN} : \OmT[[\h]] \to
(\OmT[[\h]])[1]\,,
\end{equation}
where $[1]$ denotes the shift of the total
degree by $1$\,.

Similarly, we extend the differential
$D$ (\ref{DDD}) on $\OmE[[\h]]$ to
\begin{equation}
\label{la-pi-E}
D + L_{\la_T(\pi)} : \OmE[[\h]] \to
(\OmE[[\h]])[1]\,.
\end{equation}

Notice that, the star-product $*$ in
$\OM^\h$ can
be rewritten in the form
\begin{equation}
\label{Pi}
a * b = a b +  \Pi(a,b)\,, \qquad
a, b\in \OM[[\h]]\,.
\end{equation}
where $\Pi\in \h C^1(\OM)[[\h]]$ can be
viewed as a Maurer-Cartan element of
$\Cbu(\OM)[[\h]]$\,.

Applying the map $\la_D$ (\ref{diag-D})
to $\Pi$ we get a $D$-flat Maurer-Cartan
element in $\h\,\G(M, \Cbu(\SM))[[\h]]$
and hence a new product
in $\SM[[\h]]$:
\begin{equation}
\label{diamond}
a \dia b = a b +  \la_D(\Pi)(a,b)\,, \qquad
a, b\in \G(M,\SM)[[\h]]\,.
\end{equation}
Condition \ref{uslovie} implies that
the product $\dia$ is compatible with the
following filtration on $\SM[[\h]]$
\begin{equation}
\label{filtr}
\dots \subset F^k \SM[[\h]] \subset F^{k-1} \SM[[\h]]
\subset  \dots \subset F^0 \SM[[\h]] = \SM[[\h]]\,,
\end{equation}
where the local sections of $F^k\SM[\h]$
are the following formal power series:
\begin{equation}
\label{filtr1}
\G(F^k \SM[[\h]] ) = \{\sum_{2p + m \ge k}
a_{p; i_1 \dots i_m}(x) \h^p y^{i_1} \dots y^{i_m} \}\,.
\end{equation}

Using the product $\dia$ we extend the original
differentials $D + \pa$ and $D+ \mb$
on $\OmD[[\h]]$ and $\OmC[[\h]]$ to
\begin{equation}
\label{D-dia}
D + \pa_{\dia} : \OmD[[\h]]
\to (\OmD[[\h]])[1]\,,
\end{equation}
and
\begin{equation}
\label{D-C-dia}
D + \mb_{\dia} : \OmC[[\h]]
\to (\OmC[[\h]])[1]\,,
\end{equation}
respectively.

Here $\pa_{\dia}$ (resp. $\mb_{\dia}$)
is the Hochschild coboundary (resp. boundary)
operator corresponding to (\ref{diamond}),
and $[1]$ as above denotes the shift of
the total degree by $1$\,.

Next, following the
lines of section 5.3 in \cite{thesis} we
can construct the following chain of ($\Linf$)
quasi-isomorphisms of DGLA modules:
\begin{equation}
\label{chain}
\begin{array}{ccccccc}
\Tp[[\h]] &
\stackrel{\la_T}{\,\longrightarrow\,} &
\OmT[[\h]] & \stackrel{\tcK}{\,\brarrow\,} &
\OmD[[\h]] & \stackrel{\la_D}{\,\longleftarrow\,} &
\Cbu(\OM^{\h}) \\[0.3cm]
\downarrow^{L}_{\,mod} & ~ & \downarrow^{L}_{\,mod} &
~ & \downarrow^{R}_{\,mod} & ~ & \downarrow^{R}_{\,mod} \\[0.3cm]
\cAb[[\h]] &
\stackrel{\la_{\cA}}{\,\longrightarrow\,} &
\OmE[[\h]] & \stackrel{\tcS}{\,\bblarrow\,} &
\OmC[[\h]] & \stackrel{\la_C}{\,\longleftarrow\,} &
\Cbd(\OM^\h)\,,
\end{array}
\end{equation}
where $\Tp[[\h]]$ carries the Lichnerowicz
differential
$[\pi, \,]_{SN}$, $\cAb[[\h]]$ carries
the differential $L_{\pi}$, while
$\OmT[[\h]]$, $\OmE[[\h]]$, $\OmD[[\h]]$,
$\OmC[[\h]]$ carry the differentials
(\ref{la-pi}), (\ref{la-pi-E}), (\ref{D-dia}),
(\ref{D-C-dia})\,, respectively.

The maps $\la_T$ and $\la_D$, $\la_{\cA}$,
$\la_{C}$ are
genuine morphisms of
DGLAs and their modules
as in Equations (\ref{diag-T}) and
(\ref{diag-D}). $\tcK$ is an
$\Linf$-quasi-isomorphism of the DGLAs and $\tcS$ is a
quasi-isomorphism of the corresponding
$\Linf$-modules. $\tcK$ and $\tcS$ are
obtained from $\cK$ and $\cS$
in (\ref{diag-K-Sh}), respectively, in two
steps. First, we twist\footnote{See section 2.4
in \cite{thesis} about the twisting procedure.}
$\cK$ and $\cS$ by
the Maurer-Cartan element $\la_T(\pi)$\,.
Second, we adjust them by the action of
an element $T$
of the prounipotent group
\begin{equation}
\label{mG-SM}
\mG (\,\OmD\,) = \exp [\mg(\,\OmD\,) ]
\end{equation}
corresponding to the Lie algebra
$$
\mg(\,\OmD\,) = \h\, \G(M, C^0(\SM))[[\h]] \oplus
\h\, \Om^{1}(C^{-1}(\SM)) [[\h]]\,.
$$
The element $T\in \mG (\,\OmD\,) $ is defined as
an element which transforms the
Maurer-Cartan element
\begin{equation}
\label{MC-1}
\sum_{m=1}^{\infty} \frac1{m!}
\cK_m(\la_T(\pi), \dots, \la_T(\pi))
\end{equation}
to the Maurer-Cartan element $\la_D(\Pi)$\,.

The desired trace density map (\ref{trd1})
is defined as the composition
\begin{equation}
\label{trd-def}
\trd = \Big[ H^{\bul}(\la_{\cA}) \Big]^{-1} \circ
H^{\bul}(\tcS_0)\circ
H^{\bul}(\la_C) \Big|_{HH_0(\OM^\h)}\,,
\end{equation}
where $H^{\bul}$ denotes the cohomology
functor, and $\tcS_0$ is the structure
map of the zeroth level
of the morphism $\tcS$ in (\ref{chain}).

We have to mention that
the construction of the map
(\ref{trd-def}) depends on the choice of
the element $T$ in the group (\ref{mG-SM}) which
transforms the Maurer-Cartan element (\ref{MC-1})
to $\la_D(\Pi)$, where $\Pi$ is defined in
(\ref{Pi}). Proposition \ref{abs-nonsense-mod}
implies that
altering the element $T$ changes the trace density
by the action of an automorphism
of $\OM^{\h}$ which is trivial modulo $\h$\,.
In the symplectic case all such
automorphisms are inner, while for
a general Poisson manifold there may
be non-trivial outer automorphisms.
Fortunately, we have the following proposition:
\begin{pred}
\label{nezalezhnost}
The composition (\ref{ind}) of the
trace density (\ref{trd-def}) and the
map (\ref{Ch}) is independent of the
choice of $T$ in the construction of
the trace density map.
\end{pred}
{\bf Proof.} Let $\ind$ and $\tind$
be two quantum index densities corresponding to
different choices of the element $T$ from
the group $\mG(\, \OmD \,)$ (\ref{mG-SM}).

Due to Proposition \ref{abs-nonsense-mod}
 there is an automorphism
$\tau$ of $\OM^\h$ such that
\begin{equation}
\label{1h}
\tau = 1 \qquad  m o d \qquad \h\,,
\end{equation}
and for every $\Xi \in K_0(\OM^\h)$
\begin{equation}
\label{tind}
\tind(\Xi) = \ind (\hat{\tau}(\Xi))\,,
\end{equation}
where $\hat{\tau}$ denotes the
action of $\tau$ on $K_0(\OM^\h)$\,.

But due to Proposition \ref{ind-si} the image
$\ind(\Xi)$ depends only on the principal
symbol $\si(\Xi)$\,. Therefore, since
$\tau$ does not change the principal
symbol, $\tind(\Xi) = \ind(\Xi)$\,.
 $\Box$

Let us now define the classical
index density $\ind_c$ (\ref{ind-c}) which is
a map from the K-theory of $\OM$ to
the zeroth Poisson homology of $\pi$
(\ref{pi-h}).

The well known construction of R.G. Swan
\cite{Swan} gives us the
injection\footnote{For a compact manifold $M$
this map is a bijection.}
\begin{equation}
\label{Swan-map}
\ms : K_0(\OM) \hookrightarrow K^0(M)
\end{equation}
from the $K$-theory of $\OM$ to
the $K$-theory of the manifold $M$\,.
Thus, it suffices to define the map $\ind_c$
on smooth real vector bundles.

For this, we introduce
a smooth real vector bundle $E$
over $M$ and denote
by $\End(E)$ the algebra of
endomorphisms of $E$\,.

Due to Proposition \ref{Goldman} and the formality
of the DGLA $\Cbu(\End(E))$ \cite{BGNT}, \cite{Endom}
a Maurer-Cartan element $\pi$ (\ref{pi-h}) produces a
Maurer-Cartan element $\Pi_E$ of the DGLA
$\h \, \Cbu(\End(E))[[\h]]$\,.
This element $\Pi_E$ gives us the new
associative product
\begin{equation}
\label{*-E}
a *_E \, b = a b + \Pi_E(a,b)
\end{equation}
$$
a,b \in \End(E)[[\h]]
$$
on the algebra $\End(E)[[\h]]$\,.

Due to Proposition \ref{chain-mod-twist}
the chain of quasi-isomorphisms
(\ref{diag-T}), (\ref{diag-K-Sh}),
(\ref{diag-C-E}), and (\ref{tr-cotr-tw})
for the DGLA module $(\, \Cbu(\End(E)), \Cbd(\End(E))\,)$
can be upgraded\footnote{Here we also use the fact
that every
$\Linf$-quasi-isomorphism can be replaced by a pair of
genuine (not $\Linf$) quasi-isomorphisms. }
to the chain of quasi-isomorphisms
connecting the DGLA module
$(\, \Cbu(\End(E))^{\Pi_E}, \Cbd(\End(E))^{\Pi_E}\,)$
to the DGLA module $(\Tp[[\h]], \cAb[[\h]])$ where
$\Tp[[\h]]$ carries the differential $[\pi, \, ]_{SN}$ and
$\cAb[[\h]]$ carries the differential $L_{\pi}$\,.

This chain of quasi-isomorphisms gives
us the isomorphism
$$
J_E : H_{\bul}(\Cbd(\End(E))^{\Pi_E}) \to
H_{\bul} (\cAb[[\h]], L_{\pi})
$$
from the homology of the DGLA module
$\Cbd(\End(E))^{\Pi_E}$ to the homology of
the chain complex $(\cAb, L_{\pi})$\,.

Since the complex $\Cbd(\End(E))^{\Pi_E}$
is nothing but the Hochschild chain complex for the
algebra $\End(E)[[\h]]$ with the product (\ref{*-E}),
specifying the map $J_E$ for $\bul=0$ we get
the isomorphism
\begin{equation}
\label{trd-E}
\trd_E : HH_{0} (\End(E)[[\h]], *_E) \to HP_0(M,\pi)\,,
\end{equation}
from the zeroth Hochschild homology of the
algebra $(\End(E)[[\h]], *_E)$
to the zeroth Poisson homology of $\pi$\,.

Using the map (\ref{trd-E})
we define the index density map by the equation
\begin{equation}
\label{ind-c1}
\ind_c([E]) = \trd_E ([1_E])\,,
\end{equation}
where $[1_E]$ is the class in $HH_0 (\End(E)[[\h]], *_E)$
represented by the identity endomorphism of $E$\,.

Since $\Cbu(\End(E))$ is the normalized Hochschild complex,
the group $\mG(\Cbu(\End(E)))$ acts trivially on the
$1_E$\,. Thus, Proposition \ref{abs-nonsense-mod} implies
that the map $\ind_c$ does not depend on the choices involved
in the construction of the isomorphism $\trd_E$\,.

The construction of the chain of the formality
quasi-isomorphisms (\ref{diag-T}), (\ref{diag-K-Sh}),
(\ref{diag-C-E}), and (\ref{tr-cotr-tw}) for the
DGLA module $(\, \Cbu(\End(E)), \Cbd(\End(E))\,)$
involves the choices of the connections on
the tangent bundle $TM$ and on the bundle $E$
over $M$. To show that the map $\ind_c$ (\ref{ind-c1})
is indeed well defined we need to show that
\begin{pred}
\label{nezalezh-cl}
The image $\ind_c([E])$ does not depend
on the choice of the
connections $\n$ and $\n^E$ on bundles
$TM$ and $E$\,.
\end{pred}
{\bf Proof.}
By changing the connections on $TM$ and $E$
we change the Fedosov differentials $D$ (\ref{DDD})
and $D^E$ (\ref{DDD-E1}). This means that
we twist the DGLAs  $\OmT$, $\OmD$, and
$\Omb(\Cbu(\ES))$ by Maurer-Cartan elements.
Thus, if we show that these Maurer-Cartan
elements are trivial (equivalent to zero),
the question of independence on
the connections could be easily reduced
to the application of Proposition
\ref{abs-nonsense-mod}\,.

Since the differential
(\ref{DDD}) can be viewed as a particular
case of the differential (\ref{DDD-E1})
it suffices to analyze
the differential $D^E$ (\ref{DDD-E1}).

Changing the Fedosov differential (\ref{DDD-E1})
on the DGLA $(\Omb(\Cbu(\ES)), D^E + \pa, [\,,\,]_G)$
and the DGLA
module $(\Omb(\Cbd(\ES)), D^E + \mb)$ corresponds to twisting
the DGLA structures by the Maurer-Cartan element
\begin{equation}
\label{B-E}
B^E \in \Om^1(C^0(\ES))
\end{equation}
satisfying the condition
\begin{equation}
\label{B-E-pa}
\pa B^E = 0\,.
\end{equation}

Condition (\ref{B-E-pa}) implies that
$B^E$ is a Maurer-Cartan
element of the DGLA
\begin{equation}
\label{Ker-pa}
\Om^0(\Cbu(\ES))\cap \ker \pa\, \stackrel{D^E}{\to}\,
\Om^1(\Cbu(\ES))\cap \ker \pa
 \, \stackrel{D^E}{\to} \,
 \Om^2(\Cbu(\ES))\cap \ker \pa \, \stackrel{D^E}{\to} \, \dots
\end{equation}

Thus, in virtue of Proposition \ref{Goldman},
it suffices to show that the DGLA (\ref{Ker-pa})
is acyclic in positive exterior degree.

Let $P \in \Om^{\ge 1}(\Cbu(\ES))\cap \ker \pa$
and
\begin{equation}
\label{D-E-P}
D^E P = 0\,.
\end{equation}
Let us show that an element $S\in \Omb(\Cbu(\ES))$
satisfying the equations
\begin{equation}
\label{D-E-S}
D^E S = P\,,
\end{equation}
\begin{equation}
\label{pa-S}
\pa S = 0\,.
\end{equation}
can be constructed by iterating the following
equation
\begin{equation}
\label{iter-S}
S = - \de^{-1} P +
\de^{-1}(\n S + A(S) + [\pa \ga^E, S]_G)
\end{equation}
in degrees in the fiber coordinates $y^i$\,.

Unfolding the definition of $D^E$ (\ref{DDD}),
(\ref{DDD-E1}) we rewrite the difference
\begin{equation}
\label{La}
\La = D^E S - P
\end{equation}
in the form
\begin{equation}
\label{La1}
\La = \n S - \de S + A(S) + [\pa \ga^E, S]_G - P
\end{equation}

Equation (\ref{iter-S}) implies that
$\de^{-1}S=0$ and $\chi(S)=0$\,, where
$\chi$ is defined in Equation (\ref{chi})\,.

Therefore, applying Equation (\ref{Hodge})
to $S$ we get
$$
S =\de^{-1} \de S\,.
$$
Hence,
\begin{equation}
\label{de-La}
\de^{-1}\La = 0\,.
\end{equation}

On the other hand, Equation (\ref{D-E-P})
implies that
$$
D^E \La = 0\,,
$$
which is equivalent to
\begin{equation}
\label{D-E-La}
\de \La = \n \La + A(\La) + [\pa \ga^E, \La]_{G}\,.
\end{equation}

Thus applying (\ref{Hodge}) to $\La$\,,
using Equation (\ref{de-La}) and the fact that
$\La \in \Om^{\ge 1}(\Cbu(\ES))$
we get
\begin{equation}
\label{Hodge-La}
\La = \de^{-1}(\n \La + A(\La) +
[\pa \ga^E, \La]_{G})\,.
\end{equation}
The latter equation has the unique vanishing
solution since $\de^{-1}$ raises the degree in
the fiber coordinates $y^i$\,.

The operators $\de^{-1}$ (\ref{del-1})
and $\n$ (\ref{nabla}) anticommute with $\pa$\,.
Furthermore $\pa A =0$ by definition of the form
$A$ (\ref{A}). Hence, Equation (\ref{pa-S})
follows from the definition of $S$ (\ref{iter-S}).

This concludes the proof of the proposition. $\Box$

\section{The algebraic index theorem}
Let us now formulate and prove
the main result of this paper:
\begin{teo}
\label{alg-ind}
Let $\OM^\h$ be a deformation quantization
algebra of the Poisson manifold $(M,\pi_1)$
and let $\pi$ (\ref{pi-h}) be a representative
of Kontsevich's class of $\OM^{\h}$\,.
If $\ind$ is the quantum index density
(\ref{ind}), $\ind_c$ is the classical
index density (\ref{ind-c1}) and
$\si$ is principal symbol map (\ref{prin})
then the diagram
\begin{equation}
\label{alg-ind-diag}
\begin{array}{ccccc}
K_0 (\OM^\h) &~ &\stackrel{\ind}{\longrightarrow} & ~
& HP_0(M, \h\, \pi)[[\h]]\\
~ & \searrow^{\si} &  ~  &  ~^{\ind_c}\nearrow  & ~\\
 ~ & ~ &  K_0(\OM)  & ~ & ~
\end{array}
\end{equation}
commutes.
\end{teo}
The rest of the section is devoted to
the proof of the theorem.

Let $N$ be an arbitrary natural number
and $P$ be an arbitrary
idempotent in the algebra $\Mat_N(\OM^{\h})$
of $N\times N$ matrices over $\OM^{\h}$\,.
Let $q$ be the principal symbol of $P$
$$
q = P\Big |_{\h =0}\,.
$$
Our purpose is to show that the index
$\ind([P])$ of the class $[P]$ represented
by $P$ coincides with the image $\ind_c([E])$\,,
where $E$ is the vector bundle defined by $q$\,.

Notice that $\Mat_N(\SM)$ is $\SM \otimes \End(I_N)$,
where $I_N$ denotes the trivial bundle
of rank $N$.
As in \cite{Ch-I} we would like to modify the
connection (\ref{nabla}) which is used in the
construction of the Fedosov differential
(\ref{DDD}). More precisely, we replace
$\n$ (\ref{nabla}) by
\begin{equation}
\label{nabla-q}
\n^q = \n + [\G^q,\, ] :
\SM\otimes \End(I_N) \to
\Om^1(\SM \otimes \End(I_N))\,,
\end{equation}
where
$$
\G^q = q (d q) - (d q) q\,.
$$

This connection is distinguished by the following
property:
\begin{equation}
\label{q-nabla-q}
\n^q (q) = 0\,.
\end{equation}
In general the connection $\n^q -\de + A$ is
no longer flat. To cure this problem we try
to find the flat connection within the
framework of the following ansatz:
\begin{equation}
\label{DDD-q}
D^q = D + [B^q,\, ]_{\dia}\, : \,
\SM\otimes \End(I_N)[[\h]] \to
\Om^1(\SM \otimes \End(I_N))[[\h]]\,,
\end{equation}
where $B^q\in \Om^1(\SM\otimes\End(I_N))[[\h]]$\,,
$$
B^q \Big|_{y=0} = \G^q\,,
$$
and $[\,,\,]_{\dia}$ is the commutator of
sections of $\Mat_N(\SM)[[\h]]$\,,
where the algebra $\SM[[\h]]$ is considered
with the product (\ref{diamond}).

The following proposition shows
that the desired section $B^q$
does exist:
\begin{pred}
\label{B-q}
Iterating the following equation
\begin{equation}
\label{iter-Bq}
B^q = \G^q + \de^{-1}
(\n B^q + A(B^q) + \frac1{2} [B^q,B^q]_{\dia})
\end{equation}
one gets an element
$B^q\in \Om^1(\Mat_N(\SM))[[\h]]$
satisfying the equation
\begin{equation}
\label{MC}
D B^q + \frac1{2}[B^q, B^q]_{\dia} = 0\,.
\end{equation}
\end{pred}
{\bf Proof.} First, we mention that the
process of the recursion in (\ref{iter-Bq})
converges because
\begin{equation}
\label{del-1F1}
\de^{-1} [\, F^k \Omb(\Mat_N(\SM[[\h]]))\,]
\subset F^{k+1} \Omb(\Mat_N(\SM[[\h]]))\,,
\end{equation}
where $\de^{-1}$ is defined in (\ref{del-1})
and $F^{\bul}$ is the filtration on $\SM[[\h]]$
defined in (\ref{filtr1}).

Second, since $\G^q$ does not depend on fiber
coordinates $y^i$,
$$
\de \, \G^q = 0\,.
$$
Hence, applying (\ref{Hodge}) to $\G^q$ we
get
\begin{equation}
\label{Hodge-Gq}
\G^q = \de \de^{-1}\, \G^q\,.
\end{equation}

Next, we denote by $\mu \in \Om^2(\Mat_N(\SM))[[\h]]$
the left hand side of (\ref{MC})
\begin{equation}
\label{mu-Bq}
\mu = D B^q + \frac1{2}[B^q, B^q]_{\dia}\,.
\end{equation}
Applying $\de^{-1}$ to $\mu$ we get
\begin{equation}
\label{mu-Bq1}
\de^{-1} \mu =
\de^{-1}(\n B^q + A(B^q) + \frac1{2} [B^q,B^q]_{\dia})
- \de^{-1} \de B^q\,.
\end{equation}
On the other hand, due to (\ref{nilp})
$\de^{-1} B^q = \de^{-1} \G^q$. Hence, applying
(\ref{Hodge}) to $B^q$ and using (\ref{Hodge-Gq})
we get
$$
\de^{-1} \de B^q = B^q - \G^q\,.
$$
Thus, in virtue of (\ref{iter-Bq}) and
(\ref{mu-Bq1})
\begin{equation}
\label{del-1mu}
\de^{-1} \mu  = 0\,.
\end{equation}

Using the equation $D^2=0$ it is not hard
to derive that
$$
D \mu + [B^q, \mu]_{\dia} = 0\,.
$$
In other words
$$
\de \mu = \n \mu + A(\mu) + [B^q, \mu]_{\dia}\,.
$$
Therefore, applying (\ref{Hodge}) to $\mu$
and using (\ref{del-1mu}) we get
$$
\mu = \de^{-1}(\n \mu + A(\mu) + [B^q, \mu]_{\dia})\,.
$$
The latter equation has the unique vanishing
solution due to (\ref{del-1F1}). This argument
concludes the proof of the
proposition and gives us a flat connection
of the form (\ref{DDD-q}). $\Box$

The differential $D^q$ (\ref{DDD-q}) naturally
extends to the DGLA $\Omb(\,\Cbu(\Mat_N(\SM[[\h]]))\,)$
and to the DGLA module $\Omb(\,\Cbd(\Mat_N(\SM[[\h]]))\,)$\,.
Namely,
on $\Omb(\,\Cbu(\Mat_N(\SM[[\h]]))\,)$ the differential $D^q$
is defined by the formula
\begin{equation}
\label{DDD-q-coch}
D^q = D + [\pa_{\dia} B^q,\, ]_G\,,
\end{equation}
and on $\Omb(\,\Cbd(\Mat_N(\SM[[\h]]))\,)$ it is defined by
the equation
\begin{equation}
\label{DDD-q-ch}
D^q = D + R_{\pa_{\dia} B^q }\,.
\end{equation}
Here, $\pa_{\dia}$ is the Hochschild coboundary
operator on $\Cbu(\Mat_N(\SM[[\h]]))$ where
$\SM[[\h]]$ is considered with the product
$\dia$ (\ref{diamond})\,.

Let us prove an obvious analogue of
lemma $1$ in \cite{Ch-I} (see page $10$ in
\cite{Ch-I})

\begin{lem}
\label{lemma}
If $B^q$ is obtained by iterating
Equation (\ref{iter-Bq}) then
\begin{equation}
\label{Dq0}
D\, q + [B^q, q]_{\dia} =0\,.
\end{equation}
\end{lem}
{\bf Proof.}
Since $q$ does not depend on fiber coordinates $y^i$
Equation (\ref{Dq0}) boils down to
\begin{equation}
\label{Dq01}
d\, q + [B^q, q] = 0\,,
\end{equation}
where $[\,,\,]$ stands for the ordinary
matrix commutator.

On the other hand Equation (\ref{q-nabla-q})
tells us that $d\, q = - [\G^q,q]$\,.
Thus it suffices to prove that
\begin{equation}
\label{Dq011}
[B^q, q] - [\G^q, q] = 0\,,
\end{equation}
Let us denote the right hand side of
(\ref{Dq011}) by $\Psi$
$$
\Psi = [B^q, q] - [\G^q, q]\,.
$$

Using Equations (\ref{q-nabla-q}) and (\ref{MC})
it not hard to show that
\begin{equation}
\label{D-Psi}
D \Psi + [B^q, \Psi] = 0\,.
\end{equation}
On the other hand,
Equations (\ref{nilp}) and (\ref{iter-Bq}) imply that
$\de^{-1}B^q = \de^{-1} \G^q$, and hence,
$$
\de^{-1}\Psi = [\de^{-1}B^q ,q] - [\de^{-1}\G^q, q] =0\,.
$$
Therefore, applying (\ref{Hodge}) to $\Psi$
and using (\ref{D-Psi}) we get that
$$
\Psi =
\de^{-1}
(\,\n \Psi + A(\Psi) + [B^q, \Psi]_{\dia}\,)\,.
$$
This equation has the unique vanishing
solution due to (\ref{del-1F1}).

The lemma is proved. $\Box$

We will need the following proposition:
\begin{pred}
\label{ono}
There exists an element
\begin{equation}
\label{U}
U \in Mat_N(\SM [[\h]])
\end{equation}
such that\footnote{Equation (\ref{U-mod-I}) implies
that $U$ is invertible in the algebra
$\Mat_N(\SM[[\h]])$ with the product $\dia$\,.}
\begin{equation}
\label{U-mod-I}
U = I \qquad mod \qquad Mat_N(F^1\, \SM[[\h]]\, ) \,,
\end{equation}
and
\begin{equation}
\label{U-tw}
D^q = D + [U^{-1} \dia D  U, \, ]_{\dia}\,.
\end{equation}
where $\dia$ is the obvious extension
of the product (\ref{diamond}) to $\Mat_N(\SM[[\h]])$\,.
\end{pred}
{\bf Proof.}
To prove the proposition it suffices to
construct an element $U\in Mat_N(\SM)[[\h]]$
satisfying the following equation:
$$
U^{-1} \dia D U = B^q\,,
$$
or equivalently
\begin{equation}
\label{U-Bq}
D U - U \dia B^q =0\,.
\end{equation}
We claim that a solution of (\ref{U-Bq})
can be found by iterating the equation
\begin{equation}
\label{U-iter}
U = 1 + \de^{-1}
(\n U + A(U) - U \dia B^q)\,.
\end{equation}

Indeed, let us denote by $\Phi$ the
right hand side of (\ref{U-Bq}):
$$
\Phi =  D U - U \dia B^q\,.
$$
Due to (\ref{MC})
\begin{equation}
\label{D-Phi}
D \Phi + \Phi \dia B^q =0\,.
\end{equation}
On the other hand Equations (\ref{U-iter}),
(\ref{nilp}) and (\ref{Hodge}) for $a=U$ implies that
\begin{equation}
\label{de-Phi}
\de^{-1} \Phi =0\,.
\end{equation}
Hence, applying identity (\ref{Hodge}) to
$a=\Phi$ and using (\ref{D-Phi})
we get
$$
\Phi = \de^{-1} (\n \Phi + A(\Phi) + \Phi\dia B^q)\,.
$$
Due to (\ref{del-1F1}) the latter equation has the
unique vanishing solution and the desired element
$U$ (\ref{U}) is constructed. $\Box$

Due to Equation (\ref{Dq0}) Proposition \ref{ono}
implies that the element
$$
Q = U\dia  q \dia U^{-1}
$$
is flat with respect to the initial
Fedosov differential (\ref{DDD}).
Therefore, by definition of the map
$\la_C$ (see Eq. (5.1) in chapter 5 of \cite{thesis})
\begin{equation}
\label{Q-0}
Q =\la_C(Q_0)\,,
\end{equation}
where
$$
Q_0 = Q \Big|_{y^i=0}\,.
$$
Since $Q$ is an idempotent in the
algebra $\Mat_N(\SM[[\h]])$ with the
product $\dia$ (\ref{diamond}) the element
$Q_0$ is an idempotent of the algebra
$\Mat_N(\OM^{\h})$\,.

Furthermore, due to (\ref{U-mod-I})
$$
Q_0 \Big|_{\h = 0} = q
$$
and hence, by Proposition \ref{ind-si}
$\ind([Q_0])= \ind([P])$\,.

By definition of the trace density map
$\trd$ (\ref{trd-def}) the class $\ind([Q_0])$
is represented by the cycle
\begin{equation}
\label{c}
c = \tcS_0 (\la_C(\tr\, Q_0))\,,
\end{equation}
of the complex $\OmE[[\h]]$ with
the differential (\ref{la-pi-E})\,.
Here $\tcS_0$ is the structure map of
the zeroth level of the quasi-isomorphism
$\tcS$ in (\ref{chain})\,.

Due to (\ref{Q-0})
\begin{equation}
\label{c-Q}
c = \tcS_0 (\tr\, Q)\,,
\end{equation}
where $Q = U \dia q \dia U^{-1}$\,.

Let us prove that
\begin{pred}
\label{sim}
The cycles
\begin{equation}
\label{UqU}
Q = U \dia q \dia U^{-1}
\end{equation}
and
\begin{equation}
\label{ryad}
\tQ =
\sum_{k \ge 0} (-1)^k
[q \otimes (B^q)^{\otimes(2k)} +
q \otimes (B^q)^{\otimes(2k+1)} ]
\end{equation}
are homologous in the complex
$\Omb(\,\Cbd(\Mat_N(\SM[[\h]])) \,)$
with the differential $D + \mb_{\dia}$\,.
\end{pred}
{\bf Proof.} A direct computation shows that
$$
Q - \tQ = D \psi + \mb_{\dia} \psi
$$
where
$$
\psi = \sum_{k \ge 0} (-1)^k
[ U \otimes (B^q)^{\otimes(2k)} \otimes q \dia U^{-1}
+
U  \otimes (B^q)^{\otimes(2k+1)} \otimes q \dia U^{-1} ]\,. \quad \Box
$$

It is not hard to show that
$$
\tQ = \exp(R_{B^q})\, (q)
$$

Therefore, the class $\ind([P])$
is represented by the cycle
\begin{equation}
\label{c-prime}
c' = \tcS_0 \circ \tr \circ \exp (R_{B^q}) \,(q)\,.
\end{equation}

Let us consider the following diagram of quasi-isomorphisms
of DGLA modules:
\begin{equation}
\label{tsep}
\begin{array}{ccccc}
\OmT[[\h]] & \stackrel{\tcK}{\,\brarrow\,} &
\OmD[[\h]] & \stackrel{\cotr'}{\,\longrightarrow\,} &
\Omb(\,\Cbu(\Mat_N(\SM[[\h]])) \,)  \\[0.3cm]
\downarrow^{L}_{\,mod} &
~ & \downarrow^{R}_{\,mod} & ~ & \downarrow^{R}_{\,mod} \\[0.3cm]
\OmE[[\h]] & \stackrel{\tcS}{\,\bblarrow\,} &
\OmC[[\h]] & \stackrel{\tr'}{\,\longleftarrow\,} &
\Omb(\,\Cbd(\Mat_N(\SM[[\h]])) \,)\,,
\end{array}
\end{equation}
where $\Omb(\,\Cbu(\Mat_N(\SM[[\h]])) \,)$
and $\Omb(\,\Cbd(\Mat_N(\SM[[\h]])) \,)$ carry
respectively the differentials (\ref{DDD-q-coch})
and (\ref{DDD-q-ch}), the rest DGLAs and
DGLA modules carry the same differentials
as in (\ref{chain}), and
\begin{equation}
\label{withprime}
\cotr' = \exp (-[B^q\,, \, ]_G) \circ \cotr\,,
\qquad
\tr' = \tr \circ \exp (R_{B^q})\,.
\end{equation}

Recall that $E$ is the vector bundle corresponding to
the idempotent $q$ of the algebra $\Mat_N(\OM)$\,.
In other words the rank $N$ trivial bundle
is the direct sum
\begin{equation}
\label{I-N}
I_N = E \oplus \bE\,,
\end{equation}
where $\bE$ is the bundle corresponding to $1-q$\,.

In a trivialization compatible with the
decomposition (\ref{I-N}), the endomorphism
$q$ is represented by the constant matrix
\begin{equation}
\label{q}
\tilde{q} =
\left(
\begin{array}{cc}
I_m & 0 \\
0 &  0
\end{array}
\right)\,,
\end{equation}
where $m$ is the rank of the bundle $E$\,.

If we choose different trivializations on $I_N$
the element $B^q$ in (\ref{DDD-q}) may no longer be regarded
as a one-form with values in $\SM\otimes \End(I_N)[[\h]]$\,.
Instead $B^q$ is the sum
$$
B^q = \G^q + \widetilde{B}^q
$$
of the connection form $\G^q$ (\ref{nabla-q}) and
an element $\widetilde{B}^q \in \Om^1(\SM\otimes \End(I_N))[[\h]]$
such that
$$
\widetilde{B}^q \Big|_{y=0} = 0\,.
$$
However the maps $\cotr'$ and $\tr'$ (\ref{withprime})
in Diagram (\ref{tsep}) are still well defined.
The latter follows from the fact that we deal
with normalized Hochschild (co)chains.

Lemma \ref{lemma} implies that, in a trivialization
compatible with the decomposition (\ref{I-N}),
the form $B^q$ is represented  by the block
diagonal matrix
\begin{equation}
\label{B-q-dec}
B^q =
\left(
\begin{array}{cc}
A^E & 0 \\
0 &  A^{\bE}
\end{array}
\right)\,,
\end{equation}
where $A^E$
$$
A^E = \G^E + \widetilde{A}^E\,,
$$
$$
A^{\bE} = \G^{\bE} + \widetilde{A}^{\bE}\,,
$$
$\G^E$ (resp. $\G^{\bE}$) is a connection form
of $E$ (resp. $\bE$) and
$\widetilde{A}^E \in \Om^1(\SM\otimes \End(E))[[\h]]$\,,
$\widetilde{A}^{\bE} \in \Om^1(\SM\otimes \End(\bE))[[\h]]$\,.

The latter implies that the cycle $c'$ (\ref{c-prime})
$$
c' =  \tcS_0 \circ \tr \circ \exp (R_{A^E}) \,(1_E)\,.
$$
Hence $c'$ represents the class
$H_{\bul}(\la_{\cA})(\trd_E([1_E]))$ in
the cohomology of the complex
$(\OmE[[\h]]$ with the differential $D + L_{\la_T(\pi)}$\,.
Here $\la_{\cA}$ is the embedding of $\cAb[[\h]]$
into $\OmE[[\h]]$ (\ref{diag-T})
and $H_{\bul}$ denotes the cohomology functor.

Since $c'$ is cohomologous to $c$ (\ref{c})
the statement of Theorem \ref{alg-ind} follows. $\Box$

~\\
{\bf Remark.} Theorem \ref{alg-ind} can be
easily generalized to the deformation quantization
of the algebra $\OM^{\bbC}$ of smooth complex valued
functions. In this setting we should use
the corresponding analogue of the formality
theorem from \cite{Endom} for smooth complex
vector bundles.

\section{Concluding remarks}

\subsection{The relation to the cyclic version
of the algebraic index theorem}

In \cite{TT} D. Tamarkin and B. Tsygan, inspired
by the Connes-Moscovici higher index formulas \cite{CM},
suggested the first version of the algebraic index theorem
for a Poisson manifold. This version is based on the
cyclic formality conjecture \cite{Tsygan}.

The statement of this conjecture \cite{Tsygan}
(see conjecture 3.3.2) would provide us with the cyclic
version of the trace density map which is a map
\begin{equation}
\label{trd-cyc}
\trd^{cyc}_{\bul} : HC^{per}_{\bul}(\OM^{\h})
\to H_{\bul}(\Omb(M)((u))[[\h]],\, u\, d )
\end{equation}
from the periodic cyclic homology $HC^{per}_{\bul}(\OM^{\h})$
of the deformation quantization algebra $\OM^{\h}$
to the homology of the complex
\begin{equation}
\label{complex}
(\Omb(M)((u))[[\h]],\, u\, d )
\end{equation}
where $u$ is an auxiliary variable of degree $-2$ and
$d$ is the De Rham differential.

The algebraic index theorem \cite{TT} of D. Tamarkin
and B. Tsygan describes the map (\ref{trd-cyc})
in terms of the principal symbol map
\begin{equation}
\label{symvol}
\si_{cyc}: HC^{per}_{\bul}(\OM^{\h}) \to
HC^{per}_{\bul}(\OM)
\end{equation}
and characteristic classes of $M$\,.

In order to show how our quantum index density
(\ref{ind}) fits into the picture of D. Tamarkin
and B. Tsygan let us recall that the map
(\ref{trd-cyc}) is the composition of two
isomorphisms:
$$
\trd^{cyc}_{\bul}=
\beta \circ \ttrd^{cyc}_{\bul}\,.
$$

The first isomorphism is
the map\footnote{It is the cyclic formality
conjecture which would imply the existence
of the isomorphism (\ref{trd-cyc1}).}
\begin{equation}
\label{trd-cyc1}
\ttrd^{cyc}_{\bul} :
 HC^{per}_{\bul}(\OM^{\h})
\to H_{\bul}(\Omb(M)((u))[[\h]],\, L_{\pi} + u\, d )
\end{equation}
from the periodic cyclic homology of
$\OM^{\h}$ to the homology of the complex
\begin{equation}
\label{complex1}
(\Omb(M)((u))[[\h]],\, L_{\pi} + u\, d )\,,
\end{equation}
where $L_{\pi}$ denotes the Lie derivative
along the bivector $\pi$ (\ref{pi-h}).

The second isomorphism
\begin{equation}
\label{Boris}
\beta:
  H_{\bul}(\Omb(M)((u))[[\h]],\, L_{\pi} + u\, d )
\to
  H_{\bul}(\Omb(M)((u))[[\h]],\, u\, d )
\end{equation}
is induced by the following map between
complexes (\ref{complex}) and (\ref{complex1})
$$
c \to \exp(u^{-1} \, i_{\pi})\,c\,,
$$
where $i_{\pi}$ denotes the contraction with
the bivector $\pi$\,.

The quantum index density (\ref{ind})
fits into the following commutative diagram:
\begin{equation}
\label{diag}
\begin{array}{ccc}
K_0(\OM^{\h}) & \stackrel{\ch^-_0}{\longrightarrow}&
HC^-_0(\OM^{\h})  =  HC^{per}_0(\OM^{\h})\\[0.3cm]
\downarrow^{\,\ind} &~&
\downarrow^{\,\ttrd^{cyc}_{0}}\\[0.3cm]
HP_0(M,\pi) & \stackrel{u = 0}{\longleftarrow}
& H_0(\Omb(M)[[u]][[\h]],\, L_{\pi} + u\, d )\\[0.3cm]
~&~& \downarrow^{=}\\[0.3cm]
~&~& H_0(\Omb(M)((u))[[\h]],\,
L_{\pi} + u\, d )  \\[0.3cm]
~&~& \downarrow^{\beta}\\[0.3cm]
~&~& H_0(\Omb(M)((u))[[\h]],\, u\, d )\,,
\end{array}
\end{equation}
where $\ch^-_0$ denotes the Chern character map
(see Proposition 8.3.8 in Section 8.3 in \cite{Loday})\,.

We would like to mention recent paper \cite{CF}
by A. S. Cattaneo and G. Felder. In this paper 
the authors consider a manifold $M$ equipped with a volume form
and construct an $L_{\infty}$ morphism
from the DG Lie algebra module
$CC_{\bul}^-(\OM)$ of negative cyclic chains of $\OM$
to a DG Lie algebra module modeled
on polyvector fields using the volume form. 
Although this $L_{\infty}$ morphism is not a quasi-isomorphism 
one can still use it to construct  
a specific trace on the deformation quantization algebra of a 
unimodular Poisson manifold. It would be interesting 
to find a formula for the index map corresponding 
to this trace.

\subsection{Lie algebroids and the algebraic index
theorem}
There are
two ways when the Lie algebroid theory comes in the game.
The first one, more direct, is based on the
formality theorem
for Lie algebroids (see  \cite{Cal-Dolg-Halb} and
\cite{Cal-VdB} applied
to deformation quantization of the so-called
Poisson-Lie algebroid $
F\mapsto M$ with a bracket $[\,,\,]$ and anchor
$ a: \Gamma(M,F)\mapsto \Gamma(M,TM)$.
Such algebroid carries on its fibers a "Poisson
bivector" $\pi^F\in \Gamma(M,\Lambda^2(F))$
satisfying the Jacobi identity: $[\pi^F,\pi^F] =0.$
The corresponding version of the algebraic index theorem
in this setting could be considered as a generalization of
the results of R. Nest and B. Tsygan from \cite{NT2}.

The second way concerns the natural Poisson bracket on
the \emph{dual} vector bundle $F^* \mapsto M$ of a
Lie algebroid $F\mapsto M$.

More concretely,  we will be interested in the case
when the Lie algebroid comes as the Lie algebroid $F({\cal
G})\mapsto {\cal G}_0$ associated to a Lie groupoid ${\cal
G}\rightrightarrows {\cal G}_0$.
The dual bundle $F^*(\cal G)$
carries a natural Poisson structure which is a direct
generalization of the canonical symplectic structure on
$T^*M$ and
Lie-Poisson structure on the dual space $\mg^*$ of a Lie
algebra $\mg$. (We should remark that the simplest Lie
algebroid $TM\mapsto TM,\, a=id$ is associated with
the pair groupoid ${\cal G} = M\times M$).

Following standard definitions of \cite{NWX} we will
associate
to a Lie algebroid $F\mapsto M$ the \emph{adiabatic}
Lie algebroid
$F_{\h}\mapsto M\times I$ whose total space is the
pull-back of
$F$ and the bracket $[\,,\,]_{\h} := \h[\,,\,]$.
It is interesting
and important result that this Lie algebroid comes as a Lie
algebroid of Connes's tangent groupoid $\cal G^T$ (see
\cite{Connes}).

There are two $C^*$-algebras that can be
considered in this situation. The first one
is Connes's $C^*$-algebra
$C^*(\cal G)$ of a Lie groupoid and the second one
is its "classical" counterpart --
the Poisson algebra $C_0(F^*(\cal G))$
of continuous functions on
$F^*(\cal G)\,.$ An appropriate type of
a deformation quantization
in the $C^*$-algebra context was proposed by
M. Rieffel \cite{Rief}. Let us remind
it here omitting some technical details

\begin{defi} A $C^*$-algebraic deformation
quantization (or the \emph{"strict" Rieffel's
quantization}) of a Poisson manifold $M$ is a
continuous family of
$C^*$-algebras $(A,A_{\hbar},\hbar\in I=[0,1])$
such that $A_0 = C_0(M)$ and the Poisson algebra
$\tilde A_0$ is dense in $C_0(M)$.
There is a family of sections
$$
\mF: I\mapsto
\bigsqcup_{\h\in I}A_{\h},\qquad \{
\mF(\hbar)|\mF\in A\}=A_{\hbar}
$$
and the function $\h \to ||\mF(\h)||$ continuous.
Each algebra
$A_{\hbar}$ is equipped with $*_{\hbar}$-product, a norm
$||.||_{\hbar}$ and $^{*_{\hbar}}$-involution.
The map
$$
q_h(f)=f : \tilde A_0 \mapsto A_{\hbar}
$$
satisfies the following axiom
(the "correspondence principle"):
$$
\lim_{\hbar\to 0}
||\, \frac{i}{\hbar}[\, q_{\hbar}(f), q_{\hbar}(g)\,]_{\hbar}
- q_{\hbar}( \{ f,g\})\, ||_{\hbar}=0\,.
$$
\end{defi}

The proper $C^*$ analog of the algebra $A[[\hbar]]$ is a
$C[I]C^*$-algebra (see \cite{Lands1}) and we will identify
$A_{\hbar}$ with $A/C[I,\hbar]A$ where
$C[I,\hbar]:=\{\mF\in C[I]|\mF(\hbar)=0\}\,.$

We will denote by $\pr_{\hbar}:A\mapsto A_{\hbar}$
the canonical projection and we will not distinguish
between $a\in A$ and the section
$a:\hbar\to \pr_{\hbar}(a)\,.$
Now there is a section map
$q:\tilde A_0 \mapsto A$ such that $q_{\hbar}
= \pr_{\hbar}\circ q\,.$

In concrete situation which we are interested in, the
$C^*-$deformation quantization was studied
by N. Landsman \cite{Lands}. Taking $\hbar\in I=[0,1]$,
for any Lie groupoid $\cal G$ the field
$$
A_0 := C_0(F^*({\cal G})),\qquad A_{\hbar}=C^*({\cal G}),
\qquad A = C^*(\cal G^{T})\,,
$$
(where $\cal G^{T}$ is the tangent groupoid of
$\cal G$) we obtain a $C^*$ algebraic deformation
quantization of the Lie algebroid $F^*(\cal G)\,.$

In this context, the arrow of ``symbol map'' $\sigma$ in
Diagram (\ref{ind-th})
(which in fact provides an isomorphism of the
$K$-groups) admits an ``inversion'':
\begin{equation}
\label{ind-a}
\ind_a := \sigma^{-1}: K_0(F^*({\cal G})) \to
K_0(C^*(\cal G))\,,
\end{equation}
which is called the \emph{analytic index map} (see
\cite{Car}, \cite{MontNist} and \cite{MontPier}).
This map plays a
key role in Connes's generalization of the
Atiyah-Singer index
theorem in the non-commutative geometry.
Proposition \ref{ind-si} (or more generally,
Rosenberg's theorem \cite{Ros}) gives us an
isomorphism of K-groups
$K_0(\tilde A[[\hbar]])\cong K_0(\tilde A)$.

In this setting we propose a plausible

\begin{conj}
\label{conj}
The maps $\ind$ and $\ind_c$ from
diagram (\ref{alg-ind-diag}) in
Theorem \ref{alg-ind}
are well defined
in the setting of the strict quantization.
Furthermore, the following diagram
\begin{equation}
\label{ind-3}
\begin{array}{ccccc}
K_0 (C^*({\cal G})) &~ &\stackrel{\ind}{\longrightarrow} & ~
& HP_0(F^*(\cal G), \pi)\\
~ & \nwarrow^{\ind_a} &  ~  &  ~^{\ind_c}\nearrow  & ~\\
 ~ & ~ &  K_0(F^*(\cal G))  & ~ & ~
\end{array}
\end{equation}
is commutative.
\end{conj}

Let us discuss this statement
in two important cases
(see \cite{Connes} and \cite{Lands}):
\begin{enumerate}

\item  when $\cal G$ is Connes's tangent groupoid
$${\cal G^T}
= {\cal G}_1\bigsqcup {\cal G}_2=
(M\times M\times (0,1])\bigsqcup
TM,$$ the corresponding Lie algebroid $F = TM$ and
$F^* = T^*M\,.$
We suppose that the manifold has a Riemannian metric
and denote by ${\cal K}(L^2(M))$ the algebra of
compact operators on the Hilbert space $L_2(M)$ of
square-integrable functions on $M$\,.

The strict quantization in this case coincides
with the Moyal deformation.

The associated  $C^*$-algebras in this case are:
$A_0 = C^*(TM)$ and $A_{\hbar} =
C^*(M\times M),\qquad \forall \hbar\in (0,1]$.
The first algebra is identified
(via the Fourier transform) with
$C_0(T^*M)$ and the second one is identified
with the algebra ${\cal K}(L^2(M))$.
Due to (\cite{Connes}, II.5, Prop.5.1),
we have the exact
sequence of $C^*$-algebras:
$$
0\longrightarrow C^*({\cal G}_1)
\longrightarrow C^*({\cal G^T})\stackrel{{\tilde\sigma}}
\longrightarrow C^*({\cal G}_2)\longrightarrow 0
$$
or in other terms
$$
0 \longrightarrow C_0((0,1])\otimes {\cal K}(L^2(M))
\longrightarrow C^*({\cal G^T})\stackrel{{\tilde\sigma}}
\longrightarrow C^*(C_0(T^*M))\longrightarrow 0
$$
and from the long exact sequence in $K$-theory we can
obtain the map
$$
{\tilde\sigma}_* : K_0({\cal G^T})\simeq K_0(C_0(T^*M))
=K^0(T^*M))\,.
$$

The map $\ind_a$ is nothing but
the Atiyah-Singer analytic
index
$$
\ind_a =\tr\circ \imath_* \circ({\tilde\sigma}_*)^{-1}:
K^0(T^*M)\to {\mathbb Z}
$$
with
$$\imath : M\times M\to {\cal G^T},\qquad \imath(x,y)
=(x,y,1),\qquad x,y\in M$$
and
$$
C^*({\cal G^T})\stackrel{\imath_*}\mapsto C^*(M\times M)
= {\cal K}(L^2(M))\stackrel{\tr}\equiv{\mathbb Z}\,.
$$

Conjecture \ref{conj} transforms in the following
commutative diagram:

\begin{equation}\label{ind-3-1}
\begin{array}{ccccc}
~ & ~ &  {\mathbb Z} & ~ & ~\\
~ & ^{\tr}\nearrow &  ~  &  ~\nwarrow^{\int}  & ~\\
K_0 ({\cal K}(L^2(M))) &~ &\stackrel{\ind}{\longrightarrow} & ~
&  H_{c}^{2n}(T^*M)\\
~ & \nwarrow^{\imath_* \circ(\sigma_*)^{-1}} &  ~  &  ~^{\ind_c}\nearrow  & ~\\
 ~ & ~ &  K^0(T^*M).  & ~ & ~
\end{array}
\end{equation}
Here we use the fact that the canonical Poisson
structure $\pi$ on
$T^*M$ is symplectic and hence
$$
HP_0(C_0(T^*M), \pi)\simeq H_{c}^{2n}(T^*M)
$$
where the De Rham cohomology with compact
support is used and the map from
$H_{c}^{2n}(T^*M)$ to $\bbZ$ is
given by the usual integral of
top degree forms over $M$\,.

\item If ${\cal G}_0$ is a point and
${\cal G}= G$ is a Lie
group then the associated Lie algebroid is nothing
but the Lie algebra ${\mathfrak g} =\Lie(G)$ and
the dual $F^* = {\mathfrak
g}^*$.  In this case the associated
$C^*$-algebras are $A_0 =
C^*(\mathfrak g)\simeq C_0({\mathfrak g}^*)$
(again, via the
Fourier transform) and $A_{\hbar} =
C^*(G), \quad \forall~ \hbar\in (0,1]$ is
the usual convolution algebra of $G$ defined by
a Haar measure.

The map $\sigma^{-1}$ is the composition
$$
K_0(C_0({\mathfrak g}^*))\stackrel{F^*}
\longrightarrow K_0(C^*(\mathfrak g))
\stackrel{\exp_*}
\longrightarrow K_0(C^*(G))\,,
$$
where $\exp :\mathfrak g \mapsto G $ is
 the usual exponential map and
the strict deformation quantization of the
Poisson-Lie structure
in ${\mathfrak g}^*$ was proposed by Rieffel
in \cite{Rief1}.

Diagram (\ref{ind-3}) takes the following form:

\begin{equation}\label{ind-3-2}
\begin{array}{ccccc}
K_0 (C^*(G)) &~ &\stackrel{\ind}{\longrightarrow} & ~
& HP_0({\mathfrak g}^*, \pi)\\
~ & \nwarrow^{\sigma^{-1}} &  ~  &  ~^{\ind_c}\nearrow  & ~\\
 ~ & ~ &  K_0(C_0({\mathfrak g}^*))  & ~ & ~
\end{array}
\end{equation}
\end{enumerate}

We would like to stress that our index theorem
(\ref{ind-th}) and
the conjectural "3-ind"-theorem (\ref{ind-3}) have,
in fact, the
same flavor of "index-without-index" theorems like the index
theorems in the theory related to Baum-Connes conjecture.
A deformation aspect of the Baum-Connes is discussed in
\cite{Lands}.

~\\

\noindent\textsc{Department of Mathematics,
University of California at Riverside, \\
900 Big Springs Drive,\\
Riverside, CA 92521, USA \\
\emph{E-mail address:} {\bf vald@math.ucr.edu}}

~\\
\noindent\textsc{Department of Mathematics,
University of Angers, \\
2, Lavoisier Boulevard,
Angers, 49045, CEDEX 01, France \\
and\\
Mathematical Physics Group, Theory Division\\
ITEP, 25, B.Tcheremushkinskaya, 117259, Moscow, Russia\\
\emph{E-mail address:} {\bf volodya@tonton.univ-angers.fr}}


\begin{thebibliography}{99}


\bibitem{AS} M. Atiyah and I. Singer,
The index of elliptic operators, I, III,
Ann. of Math. {\bf 87}, 2 (1968) 484-530, 546-609.


\bibitem{Bayen}  F. Bayen, M. Flato, C. Fronsdal,
A. Lichnerowicz, and D.
Sternheimer,  Deformation theory and quantization. I.
Deformations
of symplectic structures, Ann. Phys. (N.Y.), {\bf 111}
(1978) 61;\\ Deformation theory and quantization. II. Physical
applications, Ann. Phys. (N.Y.), {\bf 110} (1978) 111.



\bibitem{Ber}  F.A. Berezin,
Quantization,  Izv. Akad. Nauk. {\bf 38}
(1974) 1116-1175;\\ General concept of quantization, Commun. Math.
Phys. {\bf 40} (1975) 153-174.



\bibitem{BGNT} P. Bressler, A. Gorokhovsky, R. Nest, and
B. Tsygan, Deformations of Azumaya algebras, math.QA/0609575.

\bibitem{BNT} P. Bressler, R. Nest, and B. Tsygan,
Riemann-Roch theorems via deformation quantization.
I, Adv. Math. {\bf 167}, 1 (2002) 1-25.

\bibitem{BNT1} P. Bressler, R. Nest, and B. Tsygan,
Riemann-Roch theorems via deformation quantization.
II, Adv. Math. {\bf 167}, 1 (2002) 26-73.


\bibitem{JLB} J.-L. Brylinski,
A differential complex for Poisson manifolds,
J. Diff. Geom.  {\bf 28}, 1  (1988) 93--114.

\bibitem{Cal-Dolg-Halb} D. Calaque,
V. Dolgushev and G. Halbout,
Formality theorems for Hochschild chains in the Lie
algebroid setting, To appear in J. Reine Angew. Math;
arXive:math/0504372.

\bibitem{Cal-VdB} D. Calaque and M. Van den Bergh,
Hochschild cohomology and Atiyah classes, arXiv:0708.2725

\bibitem{CF} A. S. Cattaneo and G. Felder,
Effective Batalin--Vilkovisky theories,
equivariant configuration spaces and cyclic chains,
arXiv:0802.1706.


\bibitem{Ch-I} P. Chen and V.A. Dolgushev,
A Simple Algebraic Proof of the Algebraic Index Theorem,
Math. Res. Lett. {\bf 12}, 5-6 (2005) 655--671;
math.QA/0408210.

\bibitem{Connes} A. Connes, Non-Commutative Geometry,
{\it Academic Press, Inc. San Diego, CA}(1994)

\bibitem{CFS} A. Connes, M. Flato, and D. Sternheimer,
Closed star-products and cyclic cohomology,
Lett. Math. Phys. {\bf 24} (1992) 1-12.


\bibitem{CM} A. Connes and H. Moscovici,
Cyclic cohomology, the Novikov conjecture and
hyperbolic groups, Topology,  {\bf 29}, 3
(1990) 345--388.

\bibitem{FTHC} V.A. Dolgushev,
A Formality Theorem for Hochschild Chains,
Adv. Math.  {\bf 200}, 1  (2006) 51--101;
math.QA/0402248.

\bibitem{Endom} V.A. Dolgushev,
Formality theorem for Hochschild (co)chains
of the algebra of endomorphisms of a vector bundle,
Math. Res. Lett., {\bf 14}, 5 (2007) 757-767;
math.KT/0608112.


\bibitem{thesis} V.A. Dolgushev,
A Proof of Tsygan's formality conjecture
for an arbitrary smooth manifold, PhD thesis,
MIT; math.QA/0504420.


\bibitem{Fedosov} B.V. Fedosov,
Deformation quantization and index theory,
Akademie Verlag, Berlin, 1996.

\bibitem{FFS} B. Feigin, G. Felder, and B. Shoikhet,
Hochschild cohomology of the Weyl algebra and trace
in deformation quantization,
Duke Math. J. {\bf 127}, 3 (2005) 487--517;
math.QA/0311303.


\bibitem{GF} I.M. Gelfand and D.V. Fuchs,
Cohomology of the algebra of formal vector
fields, Izv. Akad. Nauk., Math. Ser.
{\bf 34} (1970) 322--337 (In Russian).


\bibitem{G-Kazh} I.M. Gelfand and D.A. Kazhdan,
Some problems of differential geometry and
the calculation of cohomologies of Lie algebras of
vector fields, Soviet Math. Dokl., {\bf 12},
5 (1971) 1367-1370.


\bibitem{Ger} M. Gerstenhaber, The cohomology structure of
an associative ring,  Annals of Math., {\bf 78} (1963)
267-288.


\bibitem{Gold} W. Goldman and J. Millson, The deformation
theory of representation of fundamental groups in compact
K\"ahler manifolds,  Publ. Math. I.H.E.S.,
{\bf 67} (1988) 43-96.


\bibitem{GR} S. Gutt and J. Rawnsley, Natural star
products on symplectic manifolds and quantum moment maps,
Lett. Math. Phys. {\bf 66} (2003) 123--139; math.SG/0304498.


\bibitem{K} M. Kontsevich, Deformation quantization of Poisson
manifolds, Lett. Math. Phys., {\bf 66} (2003) 157--216;
q-alg/9709040.


\bibitem{Koszul} J.L. Koszul,
Crochet de Schouten-Nijenhuis et cohomologie,
{\it Ast\'erisque}  (1985)  Numero Hors Serie, 257--271.

\bibitem{Lands1} N. Landsman, Quantization as a functor,
in "Quantization, Poisson Brackets and
beyond", Contemp. Math.,{\bf 315}, {\it AMS, Providence, RI}, (2002)
9-24.

\bibitem{Lands} N. Landsman, Deformation quantization and
the Baum-Connes conjecture, Comm. Math. Phys. {\bf 237} (2003),
87-103

\bibitem{Lich} A. Lichnerowicz,
Les vari\'et\'es de Poisson et leurs alg\`ebres de Lie
associ\'ees, J. Diff. Geom., {\bf 12}, 2  (1977)
253--300.

\bibitem{Loday} J.- L. Loday, Cyclic Homology,
Grundlehren der mathematischen Wissenschaften, 301.
Springer-Verlag, Berlin, 1992.

\bibitem{MontNist} B. Monthubert and V. Nistor,
A topological index theorem for manifolds with corners,
arXive:math/0507601.

\bibitem {MontPier} B. Monthubert and F. Pierrot,
Indice analytique et groupo\"{\i}des de Lie,
C.R. Acad. Scien. Paris, {\bf 325}, 2 (1997)
193--198.

\bibitem{NT} R. Nest and B Tsygan, Algebraic index theorem,
 Commun. Math. Phys. {\bf 172}, 2 (1995),  223-262.

\bibitem{NT2} R. Nest and B Tsygan,
Deformations of symplectic Lie algebroids,
deformations of holomorphic symplectic structures,
and index theorems,
Asian J. Math., {\bf 5}, 4, (2001), 599-636.

\bibitem{NWX} V. Nistor, A. Weinstein and P.Xu,
Pseudodifferential
operators on differential groupoids,
Pacific J. Math. {\bf 189}, 1 (1999) 117-152.

\bibitem{Q} D. Quillen,
Rational homotopy theory, Annals of Math., {\bf 90}, 2 (1969)
205--295.

\bibitem{Rief} M. Rieffel, Quantization and $C^*-$algebras,
Contemp. Math. {\bf 167} (1994) 67--97.

\bibitem{Rief1} M. Rieffel, Lie group convolution algebras
as deformation quantization of linear Poisson structures,
Amer. J. Math. {\bf 112} (1990) 657--686.

\bibitem{Ros} J. Rosenberg, Behaviour of $K-$theory under
quantization, in "Operator Algebras and Quantum Fields Theory", S.
Doplicher et al. (eds.), 404-413,{\it International Press,
Boston MA}, (1996).


\bibitem{Car} P. Rouse,
An analytic index for Lie groupoids,
arXiv:math/0612455.

\bibitem{Sh} B. Shoikhet,
A proof of the Tsygan formality conjecture for chains,
Adv. Math., {\bf 179}, 1 (2003) 7--37;
math.QA/0010321.

\bibitem{Swan} R. G. Swan, Vector Bundles and
Projective Modules, Trans. Amer. Math. Soc.,
{\bf 105} (1962) 264--277.

\bibitem{TT} D. Tamarkin and B. Tsygan,
Cyclic formality and index theorems.
Talk given at the Mosh\'e Flato Conference (2000),
Lett. Math. Phys.  {\bf 56}, 2  (2001) 85--97.


\bibitem{Tsygan} B. Tsygan,
Formality conjectures for chains,
Differential topology, infinite-dimensional Lie algebras,
and applications. 261--274,
{\it Amer. Math. Soc. Transl.} Ser. 2, 194,
Amer. Math. Soc., Providence, RI, 1999.



\bibitem{Ye} A. Yekutieli,
Mixed resolutions, simplicial sections and unipotent group
actions, Israel J. Math., {\bf 162} (2007) 1--27;
math.AG/0502206.

\end{thebibliography}
\end{document}